\documentclass[12pt]{amsart}
\usepackage{amsmath,amssymb,amsfonts,amscd}
\begin{document}

\theoremstyle{plain}
\newtheorem{thm}{Theorem}[section]
\newtheorem{lem}[thm]{Lemma}
\newtheorem{cor}[thm]{Corollary}
\newtheorem{prop}[thm]{Proposition}
\newtheorem{remark}[thm]{Remark}
\newtheorem{defn}[thm]{Definition}
\newtheorem{exam}[thm]{Example}

\newcommand{\U}{{\Upsilon}}
\newcommand{\ups}{{\upsilon}}

\numberwithin{equation}{section}
\newcommand{\mc}{\mathcal}
\newcommand{\mb}{\mathbb}
\newcommand{\surj}{\twoheadrightarrow}
\newcommand{\inj}{\hookrightarrow}
\newcommand{\red}{{\rm red}}
\newcommand{\codim}{{\rm codim}}
\newcommand{\half}{{\frac{1}{2}}}
\newcommand{\thalf}{{\textstyle \frac{1}{2}}}
\newcommand{\rank}{{\rm rank}}
\newcommand{\Pic}{{\rm Pic}}
\newcommand{\Div}{{\rm Div}}
\newcommand{\Hom}{{\rm Hom}}
\newcommand{\Id}{{\rm Id}}
\newcommand{\im}{{\rm im}}
\newcommand{\Spec}{{\rm Spec \,}}
\newcommand{\Sing}{{\rm Sing}}
\newcommand{\Char}{{\rm char}}
\newcommand{\Tr}{{\rm Tr}}
\newcommand{\Gal}{{\rm Gal}}
\newcommand{\wt}{{\rm wt}}
\newcommand{\Min}{{\rm Min \ }}
\newcommand{\Max}{{\rm Max \ }}
\newcommand{\sA}{{\mathcal A}}
\newcommand{\sB}{{\mathcal B}}
\newcommand{\sC}{{\mathcal C}}
\newcommand{\sD}{{\mathcal D}}
\newcommand{\sE}{{\mathcal E}}
\newcommand{\sF}{{\mathcal F}}
\newcommand{\sG}{{\mathcal G}}
\newcommand{\sH}{{\mathcal H}}
\newcommand{\sI}{{\mathcal I}}
\newcommand{\sJ}{{\mathcal J}}
\newcommand{\sK}{{\mathcal K}}
\newcommand{\sL}{{\mathcal L}}
\newcommand{\sM}{{\mathcal M}}
\newcommand{\sN}{{\mathcal N}}
\newcommand{\sO}{{\mathcal O}}
\newcommand{\sP}{{\mathcal P}}
\newcommand{\sQ}{{\mathcal Q}}
\newcommand{\sR}{{\mathcal R}}
\newcommand{\sS}{{\mathcal S}}
\newcommand{\sT}{{\mathcal T}}
\newcommand{\sU}{{\mathcal U}}
\newcommand{\sV}{{\mathcal V}}
\newcommand{\sW}{{\mathcal W}}
\newcommand{\sX}{{\mathcal X}}
\newcommand{\sY}{{\mathcal Y}}
\newcommand{\sZ}{{\mathcal Z}}
\newcommand{\shD}{{\widehat{\mathcal D}}}
\newcommand{\A}{{\Bbb A}}
\newcommand{\B}{{\Bbb B}}
\newcommand{\C}{{\Bbb C}}
\newcommand{\D}{{\Bbb D}}
\newcommand{\E}{{\Bbb E}}
\newcommand{\F}{{\Bbb F}}
\newcommand{\G}{{\Bbb G}}
\renewcommand{\H}{{\Bbb H}}
\newcommand{\I}{{\Bbb I}}
\newcommand{\J}{{\Bbb J}}
\newcommand{\M}{{\Bbb M}}
\newcommand{\N}{{\Bbb N}}
\renewcommand{\P}{{\Bbb P}}
\newcommand{\Q}{{\Bbb Q}}
\newcommand{\R}{{\Bbb R}}
\newcommand{\T}{{\Bbb T}}
\newcommand{\V}{{\Bbb V}}
\newcommand{\W}{{\Bbb W}}
\newcommand{\X}{{\Bbb X}}
\newcommand{\Y}{{\Bbb Y}}
\newcommand{\Z}{{\Bbb Z}}
\newcommand{\td}[1]{t_1,\dotsc,t_{#1}}
\def\ve{\,|\,}
\def\thi{t^{1/2-i}}
\def\tdn{t_1,\dots,t_n}
\newcommand{\tsdn}{t_{\sigma(1)},\dotsc,t_{\sigma(n)}}
\def\iln{i_1<\dots<i_n}

\def\sper{\sum_{\sigma\in\mathfrak{S}(n)}}
\def\P#1{\Pi_{#1}}
\def\G#1{\Gamma_{#1}}
\def\Pt#1{\Pi^\bullet_{#1}}
\def\Po#1{\Pi^\circ_{#1}}
\def\Gt#1{\Gamma^\bullet_{#1}}
\def\Go#1{\Gamma^\circ_{#1}}
\def\Ga{\Gamma}
\def\on{\{1,\dots,n\}}
\def\lpi{\ell(\pi)}
\def\g{\gamma}
\def\lga{\ell(\g)}
\def\s#1{{\Sigma}_{#1}}

\def\p#1{f^{(#1)}}
\def\pp{f'}
\def\ph#1#2#3#4{{}_2\phi_1\left[
     \begin{matrix} #1\,,\,#2\\#3 \end{matrix}\,;\,#4 \right]}

\def\ex#1#2{\left\langle\matrix #1\\ #2 \endmatrix 
\right\rangle}

\def\cbm#1#2{\left\{\begin{matrix} #1\\ #2 \end{matrix} 
\right\}_m} 

\def\eit{\ex{i}{t} }

\def\qi{(q)_\infty}
\def\i{)_\infty}
\def\la{\left\langle}
\def\ra{\right\rangle_q}

\def\l{\lambda}
\def\a{\alpha}
\def\b{\beta}
\def\T{\Theta_{11}}
\def\si{\sigma}
\def\de{\delta}
\def\Res{\operatorname{Res}}
\def\Stab{\operatorname{Stab}}

\def\tsm{\widetilde{\sum_{1\le i_1 \le \dots \le i_k \le m}}}
\def\scy{\sum_{\stackrel{\text{\rm cyclic permutations}} 
{\text{\rm of $t_1,\dots,t_k$}}}}

\def\T{\Theta}
\def\To{\T_{11}}

\def\sqb#1#2{\left[\matrix \g_1,\dots,\g_m\\{}\endmatrix
\left|\matrix #1\\ 
#2
\endmatrix \right. 
\right]}

\def\tprod{{\textstyle \prod}}
\def\tgi{{\tprod}_{j\in\ups(\g)_i}t_j}
\def\d#1{\frac{\partial}{\partial #1}}

\title{The Character of the Infinite Wedge Representation}  
\author{Spencer Bloch} 
\address{Department of Mathematics, University of Chicago, Chicago, IL 60637}
\email{bloch@math.uchicago.edu, okounkov@math.uchicago.edu}
\author{Andrei Okounkov} 
\thanks{The authors are supported in part by the NSF}
\maketitle
\setcounter{section}{-1} 

\noindent A man's character is his fate.\newline

\noindent -- Heraclitus,  {\it 540-480 B.C.}  \hfill\break

\bigskip 

\begin{abstract}
We study the character of the infinite wedge projective
representation of the algebra of differential operators
on the circle. We prove quasi-modularity of this
character and also compute certain generating functions 
for traces of differential operators which we call
correlation functions. These correlation functions are
sums of determinants built from genus 1 theta functions and 
their derivatives.
\end{abstract}  

\tableofcontents

\newpage 

\section{Introduction}\label{sec:intro} The purpose of this paper is to study the character 
associated to a certain basic projective representation, the infinite wedge representation
\cite{KP}, of the Lie algebra $\sD$ of differential operators on the circle. Concretely, we write 
\begin{equation}\label{01} D_n := \left(t\frac{d}{dt}\right)^n,\quad n\ge 0
\end{equation}
acting as differential operators on the ring $R=\C[t,t^{-1}]$ of Laurent polynomials. (Here $D_0 =
1$ is the identity operator. We occasionally write $D=D_1$.) $\sD$ has a vector space basis
$t^mD_n,\ m\in \Z,\ n\ge 0$. 

The Witt algebra of derivations of $R,\ \sW=\text{\rm Der}(R)\subset \sD$
is spanned by $t^mD,\ m\in\Z$. These algebras are graded, with $t^mD_n$ having weight $m$. One has
a central extension \cite{KR} $0\to \C\cdot c \to \shD \to \sD \to 0$ inducing the Virasoro central
extension $\widehat{\sW}$ of $\sW$, and these central extensions are also graded, with $c$ having
degree $0$. Certain graded highest weight representations of $\widehat{\sW}$ arise in connection
with conformal field theory and have been of considerable interest to physicists. These
representations have the property that for a suitable choice $\widetilde{D}\in \widehat{\sW}$
lifting $D$, the character
$$\text{\rm Trace}(q^{\widetilde{D}})
$$
is well-defined (i.e. the action of $\widetilde{D}$ is semi-simple with finite eigenspaces) and is
the $q$-expansion of a modular form \cite{ZHU}. Notice that a different choice of lifting of $D$
to $\widehat{\sW}$ multiplies the character by $q^a$ for some constant $a$. In particular, there
is at most one lifting for which the character is modular. 

For representations of $\shD$ one may consider the Abelian subalgebra
\begin{equation}\label{02} \sH := \C D_0\oplus\C D_1\oplus\C D_2\oplus\ldots\subset \sD.
\end{equation}
Suppose the action of $\sH$ is semi-simple with finite dimensional simultaneous eigenspaces. If we
choose liftings $\widetilde{D}_n \in \shD$ of $D_n$ we may define the character
\begin{equation}\label{03}
\Omega(q_0,q_1,\ldots)
:=
\text{\rm Trace}\left(q_0^{\widetilde{D}_0}q_1^{\widetilde{D}_1}q_2^{\widetilde{D}_2}\cdots
\right)\,.
\end{equation}
It is also frequently convenient to write $q_r = e^{2\pi i\tau_r}$ so the character becomes (with
 abuse of notation)
\begin{equation}\Omega(\tau_0,\tau_1,\ldots)
:= \text{\rm Trace}\left( \exp\left(2\pi i\sum_{r\ge 0}\widetilde{D}_r\tau_r\right)\right)\,.
\end{equation}
Assuming the eigenspaces of $D_1$ are themselves finite dimensional, we may specialize the $\tau_r
\mapsto 0$ for $r\ne 1$ and develop $\Omega(\tau_0,\tau_1,\ldots)$ in a formal Taylor series
expansion in
$\tau_0,\tau_2,\ldots$ with coefficients functions in $\tau_1$:
\begin{equation}\label{05} \Omega(\tau_0,\tau_1,\ldots) = \sum_A\omega_A(\tau_1)\tau^A/A!\,.
\end{equation}
Here $A=(a_0,a_2,a_3,\ldots)$ with almost all $a_j=0$, and $\tau^A/A!$ is multi-index notation. 
Finally, $\widetilde{D}_0$ is central in $\shD$ so its eigenspaces are stable under the $\shD$ action.
We write 
\begin{equation}\label{06}V(q_1,q_2,\ldots) = \text{\rm Coeff. of $q_0^0$ in }\Omega(q_0,q_1,\ldots)
\end{equation}
for the character of the $\widetilde{D}_0=0$ eigenspace. Again, by abuse of notation, we also
write $V(\tau_1,\tau_2,\ldots)$ and we expand in a series
\begin{equation}\label{07} V = \sum_{B=(b_2,b_3,\ldots)} v_B(\tau_1)\tau^B/B!\,,
\end{equation}
where, by construction, 
\begin{equation}\label{07a}
\frac{v_B(\tau_1)}{(2\pi i)^{b_2+b_3+\dots}}=
\text{\rm Trace}\big|_{\widetilde{D}_0=0} \left(
q_1^{\widetilde{D}_1} \prod_{k=2}^\infty
\left(\widetilde{D}_k\right)^{b_k} 
\right)\,.
\end{equation}
Here the trace is taken in the $\widetilde{D}_0=0$ subspace.

In the case of the infinite wedge representation, these characters have the following shape
\begin{multline}\label{08} \Omega(q_0,q_1,q_2,\ldots) = \\
q_1^{-\xi(-1)}q_3^{-\xi(-3)}\cdots 
\prod_{r\ge 0}
(1+q_0q_1^{r+\half}q_2^{(r+\half)^2}\cdots)(1+q_0^{-1}q_1^{r+\half}q_2^{-(r+\half)^2}\cdots)
\end{multline} 
\begin{multline}\label{09}V(q_1,q_2,\ldots) = 
q_1^{-\xi(-1)}q_3^{-\xi(-3)}\cdots \\
\times\sum_{(m|n)}q_1^{\sum_j (m_j+\half)+(n_j+\half)}q_2^{\sum_j (m_j+\half)^2-(n_j+\half)^2}
q_3^{\sum_j (m_j+\half)^3+(n_j+\half)^3}\cdots
\end{multline}
Here 
\begin{equation}\label{010}\xi(s) := \sum_{n\ge 1}\left(n-\thalf\right)^{-s} = (2^s-1)\zeta(s).
\end{equation}
Also, the sums are over all 
\begin{gather}\label{011}(m|n) = (m_1,m_2,\ldots,m_a|n_1,\ldots,n_a) \\
m_a>m_{a-1}>\cdots >m_1\ge 0,\quad n_a>n_{a-1}>\cdots>n_1\ge 0\notag \,.
\end{gather}
Such data are {\it Frobenius coordinates} for partitions $\lambda$ (cf. \cite{MacD}, p. 3.)  In
particular, the exponents in \eqref{09} are functions on the set of partitions. We define
\begin{equation}\label{012}p_r(\lambda) := \sum_j
  \left(m_j+\thalf\right)^r+
(-1)^{r+1}\left(n_j+\thalf\right)^r \,.
\end{equation}
Thus $p_0(\lambda)=0$ and $p_1(\lambda) = |\lambda|$ is the number being partitioned. 

The formulas \eqref{08} and \eqref{09} appear (without the anomaly terms) in \cite{FKRW} and
\cite{AFMO}. The anomaly terms depend on a choice of liftings of the $D_n$ to $\shD$. For
representations of the Virasoro algebra, this lifting is determined by modularity. The correct
analog, quasimodularity, for representations of $\shD$ was suggested by Dijkgraaf \cite{D} based
on an interpretation of $V(\tau_1,\tau_2,0,0,\ldots)$ as a generating function counting covers of a
fixed elliptic curve E with given degree and genus. From this point of view, a mirror symmetry
argument suggested an interpretation in terms of (nonholomorphic) modular forms on the dual moduli
space. More precisely, note that in \eqref{05} and \eqref{07} the variable $\tau_1$ plays a
distinguished role. The analogue of modularity involves the behavior of $\Omega$ or $V$
under the transformations
\begin{equation}\label{12} \tau_1 \mapsto
\frac{a\tau_1+b}{c\tau_1+d},\qquad \tau_j \mapsto \frac{\tau_j}{(c\tau_j+d)^{j+1}},\ j\ne 1\,,
\end{equation}
for $\begin{pmatrix}a & b \\ c & d \end{pmatrix}\in \Gamma$, where $\Gamma\subset \text{\rm SL}_2(\Z)$
is some subgroup of finite index. Given $C=(c_0,c_2,c_3,\ldots)$ (resp. $C=(c_2,c_3,\ldots)$) with
almost all $c_i=0$ define the weight of $C$
\begin{equation}\label{0013} \text{\rm wt}(C) := \sum (i+1)c_i \,.
\end{equation}
To say that the $\omega_A(\tau_1)$ (resp. the $v_B(\tau_1)$) were modular of weight
$\text{\rm wt}(A)+w$ (resp. $\text{\rm wt}(B)+w$) for some $w$ would imply that $\Omega$ (resp. $V$) was
modular of weight $w$ under the transformation \eqref{12}. This doesn't happen in the examples we
know. Instead, the $\omega_A$ and $v_B$ are {\it quasimodular} of these weights.  The notion of
quasimodular form is developed in \cite{KZ}, where a rigorous proof of quasimodularity for the
$\omega_A$ (resp. $v_B$) is given for $A=(a_0,a_2,0,0,\ldots)$ (resp. $B=(b_2,0,0,\ldots)$). By
definition, the ring of quasimodular forms is the graded algebra generated over the ring of
modular forms by the Eisenstein series of weight $2$
$$ G_2(q) := -\frac{B_2}{4}+\sum_{n=1}^\infty\sum_{d|n}dq^n.
$$
Unlike the ring of modular forms, this ring is closed under the derivation $q\frac{d}{dq}$. One
can canonically associate to a quasimodular form a certain nonholomorphic ({\it almost
holomorphic} in the terminology of op.\ cit.) modular form. 
\begin{defn}\label{defn01} A series $F = \sum_A f_A(\tau_1)\tau^A/A!$ for $A=(a_0,a_2,\ldots)$ is
quasimodular of weight $w$ if each $f_A(\tau_1)$ is quasimodular of
weight $w+\text{\rm \rm wt}(A)$. 
\end{defn}
\begin{thm}\label{th02}$\Omega(\tau_0,\tau_1,\ldots)$ (resp. $V(\tau_1,\tau_2,\ldots)$) is
quasimodular of weight $0$ (resp. weight $-\half$).
\end{thm}
Although the anomaly factors $q_{2r+1}^{-\xi(-2r-1)}$ are uniquely determined by quasimodularity,
they can also be understood solely in representation-theoretic terms. Roughly speaking, if one
tries to construct the infinite wedge representation as a representation of $\sD$ (rather than
$\shD$) one is lead to meaningless infinite constants of the form $\sum_{n\ge 1} \left(n-\half\right)^r,\
r\ge 1$. Regularizing these sums via analytic continuation of $\xi(s)$ and using the fact that
$\xi(-2r)=0$ for $r\ge 0$ leads to the stated values. 
 
The analogy between \eqref{08} and the triple product formulas for the
genus 1 theta functions continues with the following elliptic
transformation formula for $\Omega(\tau_0,\tau_1,\ldots)$.
\begin{thm}\label{th03} Define the transformation $T$ by
$$
T(\tau_j) = \tau_j -\binom{j+1}{1}\tau_{j+1}+\binom{j+2}{2}\tau_{j+2}
-\ldots\,,
$$
and set $q'_i=\exp(2\pi i T(\tau_i))$. Then we have
$$
\Omega(q_0',q_1',q_2',\ldots) =
q_0q_1^{-\half}q_2^{+\frac{1}{3}}q_3^{-\frac{1}{4}}\cdots\Omega(q_0,q_1,\ldots)
$$
\end{thm}
  
The above results constitute a pleasant but perhaps not terribly surprising generalization of the
work in \cite{D} and \cite{KZ}. The character $V$ possesses some other
hidden structure
and relation to the genus 1 theta functions,
which seems to us quite surprising and new. 

The idea is to replace
the bulky generating function \eqref{09} for the quantities
\eqref{07a} by other generating functions, where the 
auxiliary variables are attached not
to the exponents $b_i$ in \eqref{07a}
but to the indices $k$ of $\widetilde{D}_k$. These
new generating functions admit a neat evaluation in terms of
genus 1 theta functions. 

We call these generating functions $n$-point \emph{correlation
functions} because both by their definition (as averages of product of
$n$  generating series) and by their analytic structure (determinants built from
theta functions and their derivatives) these functions closely
resemble correlation functions in QFT. However, we were not able
to find any precise connection and in the present text we work with
these functions using solely the classical methods of analysis.

Concretely, the definition of these $n$-point functions is the
following. Let $f(\lambda)$ be a function on partitions, and define
\begin{equation}\label{015} \langle f\rangle_q \ := 
\sum_\lambda f(\lambda)q^{|\lambda|}\Big/ \sum_\lambda
q^{|\lambda|}\,.
\end{equation}
That is, $\langle f\rangle_q$ is the expectation of the
function $f$ provided the the probability of each partition $\l$
is proportional to $q^{|\l|}$. 
Then, the equations \eqref{07a}, \eqref{09} and \eqref{011} can be
restated as follows:
\begin{equation}\label{015a}
\text{\rm Trace}\big|_{\widetilde{D}_0=0} \left(
q_1^{\widetilde{D}_1} \prod_{k=2}^\infty
\left(\widetilde{D}_k\right)^{b_k} 
\right)=\eta(q_1) \left\langle \prod_{k=2}^\infty
\left(p_k(\l)-\xi(-k)\right)^{b_k} \right\rangle_{q_1}\,.
\end{equation}
One checks (see section \ref{sec:prepart}) that
\begin{equation}\label{016a}
p_k(\l)-\xi(-k)=\left.\left(t\frac{d}{dt}\right)^k 
\left(\sum_{i=1}^\infty t^{\l_i-i+1/2} - \frac1{\log t}\right)
\right|_{t=1}\,.
\end{equation}
By definition, set for $n=1,2,3,\dots$
\begin{equation}\label{016}
F(t_1,\dots,t_n;q)= \left\langle \prod_{k=1}^n \left(\sum_{i=1}^\infty
t_k^{\lambda_i-i+\half}\right)  \right\rangle_{q} \,, 
\end{equation}
that is, $F(t_1,\dots,t_n;q)$ 
is the expectation of the product of $n$ generating
series for the quantities \eqref{016a}. 
By linearity, all quantities
\eqref{015a} satisfying $b_2+b_3+\dots\le m$ for some $m$
are encoded in the functions \eqref{016} with $n\le m$. 
\begin{defn}\label{defn02} The functions \eqref{016} are called
$n$-point correlation functions. 
\end{defn}

The equation \eqref{016a} can be restated as follows. 
Write $t_k = e^{u_k}$ and
define a differential operator
\begin{equation}\label{017a}
\delta(u) := \frac1{u} +(2\pi i)^{-1}\sum_{r=1}^\infty
\frac{\partial}{\partial\tau_r}u^r/r!\,.
\end{equation}
Then
\begin{equation}\label{017}
\Big<\prod_{k=1}^n \Big(\sum_{i=1}^\infty
t_k^{\lambda_i-i+\half}\Big)\Big>_{q_1}
= \eta(q_1)\delta(u_1)\circ\ldots\circ\delta(u_n)V|_{\tau_2=\tau_3=\dotsc = 0}.
\end{equation}
This is also proved in section \ref{sec:prepart}. 

The form of the singular term  $1/u$  and in \eqref{017a} and
the corresponding term in \eqref{016a} very much depends on our choice of
anomaly factors in the character. 

Write $q=q_1$, write $F(t_1,\dotsc,t_n)$ for the $n$-point
function, and write $\Theta(t)=\Theta(t;q)$ for the 
following theta function: 
$$ \Theta(t) := 
\eta(q)^{-3}\sum_{n\in \Z}
(-1)^nq^{\frac{(n+\half)^2}{2}}x^{n+\half}= (q)^{-2}_\infty(x^\half -
x^{-\half})(qx)_\infty(q/x)_\infty\,.
$$
Let $\Theta^{(p)}(t) = (t\frac{d}{dt})^p\Theta(t)$. Then our main result on the $n$-point function
is following:
\begin{thm}\label{thm04}
\begin{equation}\label{018}
\vspace{-2 \jot}
F(\tdn)= 
\sum_{\sigma\in\mathfrak{S}(n)}\, 
\frac
{\displaystyle \det\left( \frac{\displaystyle \T^{(j-i+1)}(t_{\sigma(1)}\cdots
    t_{\sigma(n-j)})}{\displaystyle (j-i+1)!} 
\right)} 
{\displaystyle \T(t_{\sigma(1)})\,\T(t_{\sigma(1)} t_{\sigma(2)}) \dots
\T(t_{\sigma(1)}\cdots t_{\sigma(n)})} \,\,. 
\end{equation} 
Here $\sigma$ runs through all permutations $\mathfrak{S}(n)$
of $\{1,\dotsc,n\}$, the
matrices in the numerator have size $n\times n$, and we define 
$1/(-n)!=0$ if $n\ge 1$. 
\end{thm}

For the $1$-point function, this is simply
$$\Big<\sum_{i=1}^\infty
t^{\lambda_i-i+\half}\Big>_q = \frac{1}{\Theta(t)}\,, 
$$
because
$$
\T'(1)=1\,.
$$
For $n=3$ the equation \eqref{018} becomes:
\begin{multline}
F(t_1,t_2,t_3)=\\
\frac1{\displaystyle\T(t_1 t_2 t_3)}
\sum_{\sigma\in\mathfrak{S}(3)}
\det
\left(
 \begin{array}{rrc}
  \frac{\displaystyle\T'\left(t_{\sigma(1)} t_{\sigma(2)}\right)}
      {\displaystyle \T\phantom{{}'}\left(t_{\sigma(1)} t_{\sigma(2)}\right)}&
  \frac{\displaystyle 1}{\displaystyle 2} \,
  \frac{\displaystyle\T''\left(t_{\sigma(1)}\right)}
    {\displaystyle\T\phantom{{}''}\left(t_{\sigma(1)}\right)}&
  \frac{\displaystyle\T'''(1)}{\displaystyle 3!}\\
  {\displaystyle 1}&
  \frac{\displaystyle\T'\phantom{{}'}\left(t_{\sigma(1)}\right)}
  {\displaystyle\T\phantom{{}''}\left(t_{\sigma(1)}\right)}& 
  0\\
  &
  {\displaystyle 1}&
  1
 \end{array}
\right)
\end{multline}
The traces \eqref{015a} can be computed from \eqref{018} 
using repeated differentiation, L'Hospital's rule and 
formulas for the derivatives $\Theta^{(p)}(1;q)$ in terms of
the Eisenstein series:
\begin{align*}
\T^{(3)}(1;q)&= - 6 G_2(q)\,,\\ 
\T^{(5)}(1;q)&= - 10 G_4(q)+ 60 G_2(q)^2\,, \\
\T^{(7)}(1;q)&= - 14 G_6(q)+420 G_4(q) G_2(q)- 840 G_2(q)^3\,, 
\end{align*}
and so on, see \eqref{th'}. Note that all even order derivatives
of the odd function $\T(x)$ vanish at $x=1$.

Sections \ref{sec:results}-\ref{sec:pf} are devoted 
to the proof of \eqref{018}.
An essential ingredient of this proof are the following $q$-difference
equations for the functions  $F(t_1,\dots,t_n)$ 
which is established in section 8:
\begin{thm}\label{thm05} For $n=1,2,\dots$ we have 
\begin{multline*}
F(qt_1,t_2,\dots,t_n)=-q^{1/2} t_1 \dots t_n \times \\ \sum_{s=0}^{n-1} 
(-1)^s 
\sum_{1<i_1<\dots < i_s \le n} F(t_1 t_{i_1} t_{i_2} \cdots t_{i_s},
\dots,\widehat{\,t_{i_1}}, \dots,\widehat{\,t_{i_s}}, \dots) \,.
\end{multline*}
\end{thm}
These equations are analogs of elliptic transformations for the 
$n$-point functions $F(t_1,\dots,t_n)$\,.

Here is another point of view about these things. Let $\sA$ be an associative algebra with $1$, and let
$\rho : \sA \to \text{\rm End}(F)$ be a representation of the associated Lie algebra (with commutator
bracket). We do {\it not} assume $\rho$ compatible with the associative algebra structures. Let
$D\in \sA$ be given and assume $\rho(D)$ is semisimple with finite dimensional eigenspaces. Let
$F_a \subset F$ be the eigenspace of $\rho(D)$ with eigenvalue $a$. Since $[D, D^n]=0$, $\rho(D^n)$
stabilizes $F_a$, and it makes sense to consider for a variable $u$ the expression
\begin{equation}\label{019}f_\rho(u,z):= \chi_\rho(z)^{-1}\sum_a\Big(\text{\rm Tr}_{F_a}\sum_{n\ge 0}
\rho(D^n)u^n/n!\Big)e^{az}\,.
\end{equation}
where $\chi_\rho(x) := \text{\rm Tr}e^{x\rho(D)}$ is the trace of $\rho$. If in fact
$\rho(D^n)=\rho(D)^n$ for all $n$, we find 
\begin{equation}f_\rho(u,z) = \chi_\rho(z)^{-1}\sum_a
\text{\rm Tr}_{F_a}e^{at}e^{az}=\chi_\rho(u+z)/\chi_\rho(z) \,. 
\end{equation}
More generally, letting
$n=(n_1,\ldots,n_r)$ run through $r$-tuples of non-negative integers and replacing $\rho(D^n)$ with
$\prod_i \rho(D^{n_i})$ and $u^n/n!$ with $\prod u_i^{n_i}/n_i!$ we may define
$f_\rho(u_1,\ldots,u_r;z)$. When $\rho$ is compatible with associative multiplication, we have
$f_\rho = \chi_\rho(u_1+\ldots+u_r+z)/\chi_\rho(z)$. 

When $\rho$ is a projective representation, $\rho(D^n)$ is well defined only up to adding
$\alpha_n\cdot \text{\rm Id}$. Such a modification replaces $f_\rho(u,z)$ with
$f_\rho(u,z)+(\sum_n\alpha_nu^n/n!)$. If $\sH$ is spanned by the powers of $D$ as in
\eqref{02} we find writing $\chi_\rho(\tau_0,\tau_1,\tau_2,\ldots)$ for the full character that
\begin{multline}\label{020} f_\rho(u_1,\ldots,u_r;\tau_1) = \\
\chi_\rho(\tau_1)^{-1}\sum_{n_1,\ldots,n_r\ge
0}\left(\prod_{i=1}^r\frac{\partial}{\partial\tau_{n_i}}\right)
\chi_\rho(\tau_0,\tau_1,\tau_2,\ldots)
|_{\tau_0=\tau_2=\ldots=0}\prod_i u_i^{n_i}/n_i! \,. 
\end{multline}
For the representation \eqref{06}, this is \eqref{017}. 

Finally, section \ref{sec:example} proves quasimodularity and computes the $n$-point function for
a representation of a subalgebra of $\sD$ studied in \cite{B}. 

The first author would like to
acknowledge considerable inspiration from conversations and correspondence with V. Kac and D.
Zagier. It was Zagier who suggested that these characters might be quasimodular.

\section{The Infinite Wedge Representation}\label{sec:rep}

Material in this section is quite well known. References are \cite{K}, \cite{K2}, \cite{PS},
\cite{FKRW}. We work throughout with vector spaces, Lie algebras, and representations of Lie
algebras over the  field
$\C$.  Recall that a (projective) representation of a Lie algebra $\sL$ is a Lie algebra
homomorphism
$\sL \to \text{\rm End}(F)$; (resp. $\sL \to \text{\rm End}(F)/\C\cdot\text{\rm Id}$) for a vector space $F$.
A projective representation gives rise to a central extension of Lie algebras via pullback:
\begin{equation}\begin{CD}0 @>>> \C\cdot c @>>> \widetilde{\sL} @>>> \sL @>>> 0 \\
@. @VVV @VVV @VVV @. \\
0 @>>>\C\cdot\text{\rm Id} @>>> \text{\rm End}(F) @>>> \text{\rm End}(F)/\C\cdot\text{\rm Id} @>>> 0\,.
\end{CD}
\end{equation}

Let $V$
be a vector space with basis
$v_n,\ n\in
\Z$. Define the finite matrices
\begin{equation}\sA_0 := \bigoplus_{i,j\in \Z} \C E_{ij}\subset \text{\rm End}(V)\,.
\end{equation}
where $E_{ij}(v_i)=v_j$ and $E_{ij}(v_k)=0;\ k\ne i$. Let $\sA\supset\sA_0$ be the larger space of
all matrices supported in a bounded strip. An infinite matrix
$M=\sum_{\text{\rm infinite}}a_{ij}E_{ij}$ lies in $A$ if and only if there exists a constant $c$ such
that $a_{ij}=0$ if $|i-j|>c$. For example, $\Id\in\sA-\sA_0$. Suppose given a representation
$\alpha_0$ of $\sA_0$ which extends to a projective representation $\alpha$ of $\sA$:
\begin{equation}\label{13} \begin{array}{ccc} \sA_0 & \stackrel{\iota}{\hookrightarrow} &
\makebox[.8in][l]{$\sA$}
\\
\rule[-4mm]{0cm}{9mm}\makebox[.1in][l]{$\downarrow\! \alpha_0$} &
\makebox[.1in]{$\swarrow\!\sigma$} &
\makebox[.7in][l]{$\downarrow\!
\alpha$}
\\
\text{\rm End}(F) & \twoheadrightarrow & \text{\rm End}(F)/\C\cdot\Id\,.
\end{array}
\end{equation}
Let $\sigma$ be a lifting of $\alpha$ which is a map of vector spaces, not necessarily compatible
with Lie algebra structures and not necessarily extending $\alpha_0$. Assume, however, that
$\sigma$ is continuous in the sense that for any $f\in F$ and any
$x=\sum_{\text{\rm infinite}}x_{ij}\iota(E_{ij})\in \sA$ there exists $N=N(x,f)$ such that
$x_{ij}\sigma\iota(E_{ij})(f)=0$ for $|i|+|j|>N$, and  
$$\sigma(x)(f) = \sum_{|i|+|j|\le N}x_{ij}\sigma\iota(E_{ij})(f)\,.
$$
Notice that such a $\sigma$ is unique up to a finite modification. More precisely, if $\sigma'$ is
another map with the same property,  $\sigma\iota(E_{ij})=\sigma'\iota(E_{ij})$ for almost all
$i,j$.  We obtain in this way a map
\begin{equation}\label{14} \epsilon := \alpha_0-\sigma\circ\iota : \sA_0 \to \C\cdot\Id\,.
\end{equation}
The functional $\epsilon$ is uniquely determined (up to a finite modification) by $\alpha_0,\alpha$. 
\begin{defn}\label{reg} Let $\sS$ be a collection of sequences of complex numbers
$c=\{c_{ij}\}_{i,j\in
\Z}$. A regularization scheme for $\sS$ is a map $r:\sS \to \C$ such that if $\{c_{ij}\},
\{c_{ij}'\}\in\sS$ differ for only finitely many pairs $i,j$, then
$$r(\{c_{ij}\}) = r(\{c_{ij}'\})+\sum_{i,j}(c_{ij}-c_{ij}').
$$
(We will usually assume $\sS$ is saturated in the sense that if $c\in \sS$ and $c'$ differs from $c$
for only finitely many $i,j$ then $c'\in \sS$.) 
\end{defn}

Assume now we have $\alpha, \alpha_0$, and that there exists a $\sigma$ as in \eqref{13} above.
Let $\widehat{\sA}$ be the corresponding central extension of $\sA$. Let $\sS,r$ be as in
definition (\ref{reg}). Suppose given $x=\sum x_{ij}\iota(E_{ij})\in \sA$, such that
$\{x_{ij}\epsilon(E_{ij})\}\in \sS$, where $\epsilon$ is defined in \eqref{14}. Then
$$\sigma(x) + r(\{x_{ij}\epsilon(E_{ij})\})\cdot\Id \in \text{\rm End}(F)
$$
is well-defined independent of the choice of $\sigma$. Thus we have a lifting $\hat x\in
\widehat{\sA}$ depending on the regularization scheme $r$ but not on the choice of lifting
$\sigma$. 
\begin{exam}\label{ex12} For the infinite wedge representation, 
\begin{equation}\epsilon(E_{ij}) = \begin{cases} 1 & i=j\le 0\,, \\ 0
  & \text{\rm else\,.}
\end{cases}
\end{equation}
The elements $x$ we want to lift have the form $x=\sum_{m\in\Z}\left(m-\thalf\right)^kE_{mm}$ for
$k=0,1,2,\ldots$, so the sequences we need to regularize are 
$$\left\{\left(m-\thalf\right)^k\right\}_{m\le 0}\,.
$$
The natural way to do this is by analytic continuation of the function 
\begin{equation}\label{16}\xi(s) = \sum_{n\ge 1} \left(n-\thalf\right)^{-s}
\end{equation}
We define $r(\{(m-\half)^k\}_{m\le 0})=(-1)^k\xi(-k)$\,. 
\end{exam}

The infinite wedge representation $F:=\Lambda^\infty V$ is the vector space with basis
\begin{equation}v_I = v_{i_m}\wedge v_{i_{m-1}}\wedge v_{i_{m-2}}\wedge\ldots\,,
\end{equation}
where $i_{k-1}<i_k$ and $i_k = k$ for $k\ll 0$. By definition $m$ is the {\it charge} of $v_I$ and 
$$\text{\rm energy}(v_I) := \sum_{k=m}^{k=-\infty}(i_k-k)\,. 
$$
For example, the elements
\begin{equation} |m\rangle := v_m\wedge v_{m-1}\wedge v_{m-2}\wedge\ldots
\end{equation}
have charge $m$ and energy $0$. The wedge here is formal, but we can use the expected alternating
linearity to make sense for $A\in \sA_0$ of expressions like
\begin{equation}\label{19} Av_I := \sum_{k=m}^{k=-\infty}v_{i_m}\wedge\ldots\wedge Av_{i_k}\wedge
v_{i_{k-1}}\wedge \ldots \quad .
\end{equation}
We obtain a representation $\alpha_0 : \sA_0 \to \text{\rm End}(F)$. To extend $\alpha_0$ to $\sA$, we
would have to define expressions like
\begin{equation}\label{110} \sum_{p\in \Z}a_pE_{p,p+r}(v_I)\,.  
\end{equation}
When $r\ne 0$ it is easy to see that $E_{p,p+r}v_I = 0$ for $|p|>>0$, so \eqref{110} makes sense.
Similarly $E_{p,p}v_I=0$ for $p>>0$. Define $\sigma : \sA_0 \to \text{\rm End}(F)$ by
\begin{equation} \sigma(E_{ij}) = \begin{cases}\alpha_0(E_{ij}) & i\ne j \text{\rm  or
}i=j>0 \\
\alpha_0(E_{ii})-\Id & i\le 0\end{cases}
\end{equation}
It is straightforward to check that for any $r$ and any $I$ there exists an $N$ such
that $\sigma(E_{p,p+r})(v_I)=0$ for $|p|\ge N$, so $\sigma$ extends to a map $\sA\to\text{\rm End}(F)$
as in \eqref{13}. Define $\alpha$ to be the composition 
$$\sA\stackrel{\sigma}{\to}\text{\rm End}(F)\to
\text{\rm End}(F)/\C\cdot\Id\,.
$$ 
\begin{lem} The map $\alpha$ is a projective representation of $\sA$. 
\end{lem}
\begin{proof} The assertion is there exists a bilinear map $a: \sA\times \sA \to \C$ with
$$[\sigma(x),\sigma(y)] = \sigma([x,y]) + a(x,y)\Id.
$$
One defines
$$a(E_{ij},E_{k\ell}) = \begin{cases}1 & i=\ell\le 0 \text{\rm\ and }j=k\ge 1\,, \\
-1 & i=\ell\ge 1 \text{\rm\ and } j=k\le 0 \,,\\
0 & \text{\rm else}\,.\end{cases}
$$
Because of the conditions on the signs of the indices, $a$ extends to $\sA\times\sA$ as desired.
\end{proof}

Given $u\in \C$, we have an action as differential operators of $\sD$ on $\C[t,t^{-1}]t^u$.
Identifying this space with $V$ via $t^{n+u} \mapsto v_n$, we get a mapping of associative
algebras (and hence a fortiori of Lie algebras)
\begin{gather*} \delta_u : \sD \to \sA\subset \text{\rm End}(V)\,, \\
\delta_u(t^pD^r) = \sum_{n\in \Z}(n+u)^rE_{n,n+p} \notag\,.
\end{gather*}
We define a projective representation
\begin{gather*}\rho_u := \alpha\circ\delta_u : \sD \to \text{\rm End}(F)/\C\cdot\Id\,, \\
\rho := \rho_{-\half}\notag\,.
\end{gather*}
We will customarily view $\rho$ as a representation on a central extension, $\rho : \shD \to
\text{\rm End}(F)$. Our objective now is to compute the character of the representation $\rho_u$.
Define operators $\psi^*_{-r},\ r\ge 0$ and $\psi_{-r},\ r>0$ on $F$ by
$$\psi_{-r}(v_I) := v_r\wedge v_I;\quad \psi^*_{-r}(\ldots\wedge v_{-r}\wedge\ldots) =
\ldots\wedge\widehat{\rule{0cm}{3mm}v_{-r}}\wedge\ldots\quad. 
$$
A basis for $F$ can be written
\begin{multline*} \psi_{-i_1}\cdots\psi_{-i_a}\psi^*_{-j_1}\cdots\psi^*_{-j_b}|0\rangle\,,\\
0<i_1<i_2<\ldots<i_a,\ 0\le j_1<j_2<\ldots<j_b\,.
\end{multline*}
We have
\begin{multline*}\sigma\circ\delta_u(D_r)\Big(\psi_{-i_1}\cdots\psi_{-i_a}\psi^*_{-j_1}
\cdots\psi^*_{-j_b}|0\rangle\Big) = \\
\sigma\begin{pmatrix}\ddots &&\\ &(n+u)^r &\\ &&\ddots\end{pmatrix}\Big(\psi_{-i_1}\cdots
 \psi_{-i_a}\psi^*_{-j_1}
\cdots\psi^*_{-j_b}|0\rangle\Big) = \\
\Big(\sum_{\ell=1}^a (i_\ell+u)^r - \sum_{m=1}^b (-j_m+u)^r\Big)\Big(\psi_{-i_1}\cdots
\psi_{-i_a}\psi^*_{-j_1}
\cdots\psi^*_{-j_b}|0\rangle\Big).
\end{multline*}
For example, writing $q_r = \exp(2\pi i\tau_r)$ we get
\begin{gather*}\exp\Big(2\pi i\sum_{r\ge 0}\tau_r\cdot\sigma\circ\delta_u(D_r)\Big)(\psi_{-n}|0\rangle)
= q_0q_1^{n+u}q_2^{(n+u)^2}\cdots \,,\\
\exp\Big(2\pi i\sum_{r\ge 0}\tau_r\cdot\sigma\circ\delta_u(D_r)\Big)(\psi^*_{-n}|0\rangle)
= q_0^{-1}q_1^{n+u}q_2^{-(n+u)^2}\cdots \,.\notag
\end{gather*}
On all of $F$ we find
\begin{multline*}\text{\rm Tr}\exp\Big(2\pi i\sum_{r\ge 0}\tau_r\cdot\sigma\circ\delta_u(D_r)\Big) = \\
\prod_{n\ge
0}\big(1+q_0q_1^{n+u+1}q_2^{(n+u+1)^2}\cdots\big)\big(1+q_0^{-1}q_1^{n-u}q_2^{-(n-u)^2}\cdots\big).
\end{multline*}
We now specialize to the case $u=-\half$. The reason why this particular value of $u$ yields a
quasimodular character is not clear, but it may have to do with the fact that for $u=p+\half\in
\Z+\half$ the space $V=\C[t,t^{-1}]t^u$ has a nondegenerate, symmetric bilinear form
$$(t^{n+u},t^{m+u}) := \text{\rm res}_{t=0} t^{n+m+2u-1}dt.
$$
Further, if we consider the involution 
$$\sigma : \sD \to \sD;\quad \sigma(t^aD^b) := -t^a(-D-a)^b,
$$
we find
$$\sigma^2 = I;\quad \sigma[x,y] = [\sigma(x),\sigma(y)];\ (x(a),b)+(a,\sigma(x)(b)) = 0.
$$

With reference
to example (\ref{ex12}), we note that
$\xi(s) = (2^s-1)\zeta(s)$ vanishes at $s=0,-2,-4,\ldots$ We therefore define the character of the
infinite wedge representation of
$\sD$ to be
\begin{multline}\label{1_18}\Omega(q_0,q_1,\ldots) = \Omega(\tau_0,\tau_1,\ldots) = \\
q_1^{-\xi(-1)}q_3^{-\xi(-3)}\cdots\prod_{n\ge
0}\big(1+q_0q_1^{n+\half}q_2^{(n+\half)^2}\cdots\big)\big(1+q_0^{-1}q_1^{n+\half}
q_2^{-(n+\half)^2}\cdots\big)\,.
\end{multline}

Note that $\delta_u(D_0)=\sum_{i\in \Z} E_{ii}$ is the charge operator
\begin{multline*}\delta_u(D_0)\psi_{-i_1}\cdots
 \psi_{-i_a}\psi^*_{-j_1}
\cdots\psi^*_{-j_b}|0\rangle = \\
(a-b)\psi_{-i_1}\cdots\psi_{-i_a}\psi^*_{-j_1}\cdots\psi^*_{-j_b}|0\rangle.
\end{multline*}
Operators in $\sD$ preserve the charge, so the charge eigenspaces are stable under $\rho_u$. For
$u=-\half$, the characters of the subrepresentation of charge $n$ is the coefficient of $q_0^n$
in \eqref{1_18}. Of particular interest is the charge $0$ part. Using \eqref{09} and \eqref{012},
this can be written
\begin{multline}\label{120} V(q_1,q_2,\ldots) = V(\tau_1,\tau_2,\ldots) := \\
q_1^{-\xi(-1)}q_3^{-\xi(-3)}\cdots\{\text{\rm Coeff. of $q_0^0$ in $\Omega(q_0,q_1,\ldots)$}\}= \\
q_1^{-\xi(-1)}q_3^{-\xi(-3)}\cdots\sum_\lambda q_1^{p_1(\lambda)}q_2^{p_2(\lambda)}\cdots\,.
\end{multline}
Here the sum is over all partitions $\lambda$.

\section{Elliptic transformation of the character $\Omega$}\label{sec:transf}

In this section we want to expose some interesting analogies between $\Omega(q_0,q_1,q_2,\ldots)$ 
\eqref{1_18} and the classical genus $1$ theta function, which we will write
\begin{equation}\theta(q_0,q_1) := \sum_{n\in\Z}q_0^nq_1^{n^2/2}.
\end{equation}
The triple product formula (cf. \cite{MUM}, p. 70) implies
\begin{multline}\label{22}\theta(q_0,q_1) = \prod_{m\ge 1}(1-q_1^m)\prod_{n\ge
0}(1+q_0q_1^{n+\half}) (1+q_0^{-1}q_1^{n+\half}) = \\
\eta(q_1)\Omega(q_0,q_1,1,1,\ldots).
\end{multline}
Here $\eta(q_1)=q_1^{\frac{1}{24}}\prod_{m\ge 1}(1-q_1^m)$ is the classical eta function, and we
have used the fact that the first anomaly factor in $\Omega$ is given by
$$q_1^{-\xi(-1)}=q_1^{-(2^{-1}-1)\zeta(-1)} = q_1^{-\frac{1}{24}}.
$$
We want to generalize to $\Omega$ the elliptic transformation law
\begin{equation}\theta(q_0q_1^{-1},q_1) = q_0q_1^{-\half}\theta(q_0,q_1).
\end{equation}
For this, we define infinite matrices
\begin{equation} U:= \begin{pmatrix}0 & -1 & 0 & 0 & \hdots \\ 0 & 0 & -2 & 0 & \hdots \\
0 & 0 & 0 & -3 & \hdots \\
\vdots & \vdots & \vdots & \vdots
\end{pmatrix}
\end{equation}
and 
\begin{equation} T:= \exp(U) = \begin{pmatrix}1 & -1 & 1 & -1 & 1 &\hdots \\\vspace{4pt}
0 & 1 & -\binom{2}{1} & \binom{3}{2} & -\binom{4}{3} & \hdots \\
0 & 0 & 1 & -\binom{3}{1} & \binom{4}{2} & \hdots \\
\vdots & \vdots & \vdots & \vdots & \vdots
\end{pmatrix} \quad. 
\end{equation}
Define for $j\ge 0$
\begin{gather*}\tau_j' := T(\tau_j) = \tau_j -\binom{j+1}{1}\tau_{j+1}+\binom{j+2}{2}\tau_{j+2}
-\ldots\,, \\
q_j' := T(q_j) = q_jq_{j+1}^{-\binom{j+1}{1}}q_{j+2}^{\binom{j+2}{2}}\cdots\,. \notag
\end{gather*}
\begin{thm}\label{thm21} $\Omega(q_0',q_1',q_2',\ldots) =
q_0q_1^{-\half}q_2^{+\frac{1}{3}}q_3^{-\frac{1}{4}}\cdots\Omega(q_0,q_1,\ldots)$.
\end{thm}
\begin{lem}\label{lem22} For $n\ge 2$ we have
$$\sum_{i=1}^{n-1}(-1)^i\binom{n}{i}\xi(i-n) = \frac{(-1)^{n+1}}{n+1}+\frac{(-1)^n}{2^n}\,.
$$
\end{lem}
\begin{proof}[Proof of lemma] Recall one has Bernoulli numbers $B_n$, $n\ge 0$ satisfying $B_0=1,\
B_1=-\half,\ \zeta(1-n)=-\frac{B_n}{n}$ for $n\ge 2$. Substituting and using the identity
$\frac{n+1}{j+1}\binom{n}{j}= \binom{n+1}{j+1}$, the desired formula becomes
\begin{equation}\label{27}\sum_{k=1}^n(-1)^{k-1}\binom{n+1}{k}(1-2^{1-k})B_k \stackrel{?}{=} 
\frac{n+1}{2^n} - 1 \,. 
\end{equation}
Consider the Bernoulli polynomials
$$B_N(x) := \sum_{k=0}^N\binom{N}{k}B_kx^{N-k} \,,
$$
which may be defined by the generating function
$$\sum_{N=0}^\infty B_N(x)t^N/N! = \frac{te^{xt}}{e^t-1}\,.
$$
As a consequence we get
\begin{equation}\label{28}B_N(x+1)-B_N(x) = Nx^{N-1} \,. 
\end{equation}
Summing for $-1\ge x\ge -p$  we get for $N=n+1$ and $p\ge 1$ the equivalent identities
\begin{gather}(n+1)\sum_{\ell=-1}^{-p}\ell^n+B_{n+1}(-p)-B_{n+1}(0)=0 \,,\notag \\
(n+1)\sum_{\ell=-1}^{-p}\ell^n+(-p)^{n+1}+\sum_{k=1}^n B_k
\binom{n+1}{k}(-p)^{n+1-k} = 0\,.
\label{29}
\end{gather}
(We are grateful to V. Kac for suggesting \eqref{29}.) 
It is straightforward to deduce \eqref{27}
from \eqref{29} taking $p=1,2$\,.
\end{proof}
\begin{proof}[Proof of theorem] The transformation $T$ satisfies 
\begin{gather*}T(q_0)T(q_1)^sT(q_2)^{s^2}\cdots = q_0q_1^{s-1}q_2^{(s-1)^2}\cdots\,, \\
T(q_0)^{-1}T(q_1)^sT(q_2)^{-s^2}\cdots = q_0^{-1}q_1^{s+1}q_2^{-(s+1)^2}\cdots\,. \notag
\end{gather*}
Note the second identity follows from the first, replacing $s$ by $-s$ and inverting. The first
identity is left for the reader. It follows that
\begin{multline}\label{211}\frac{\prod_{n\ge
0}\big(1+T(q_0)T(q_1)^{n+\half}\cdots\big)\big(1+T(q_0)^{-1}T(q_1)^{(n+\half)}
T(q_2)^{-(n+\half)^2}\cdots\big)}{\prod_{n\ge
0}\big(1+q_0q_1^{n+\half}\cdots\big)\big(1+q_0^{-1}q_1^{(n+\half)}
q_2^{-(n+\half)^2}\cdots\big)} \\
=
\frac{1+q_0q_1^{-\half}q_2^{\frac{1}{4}}q_3^{-\frac{1}{8}}
\cdots}{1+q_0^{-1}q_1^{\half}q_2^{-\frac{1}{4}} \cdots} = q_0q_1^{-\half}q_2^{\frac{1}{4}}\cdots\,.
\end{multline}
As a consequence of the lemma we have
\begin{equation}\label{212}
\frac{T(q_1)^{-\xi(-1)}T(q_3)^{-\xi(-3)}\cdots}{q_1^{-\xi(-1)}q_3^{-\xi(-3)}\cdots} =
q_2^{\frac{1}{3}-\frac{1}{4}}q_3^{-\frac{1}{4}+\frac{1}{8}}q_4^{\frac{1}{5}-\frac{1}{16}}
\cdots \,.
\end{equation}
The proof follows by combining \eqref{211} and \eqref{212}. 
\end{proof}

Recall we have defined $V(q_1,q_2,\ldots)$ to be the coefficient of $q_0^0$ in
$\Omega(q_0,q_1,\ldots)$. 
\begin{thm}\label{thm23} We have the following series expansion for $\Omega$:
$$\Omega(q_0,q_1,q_2,\ldots) = \sum_{n=-\infty}^{n=\infty}V(T^{-n}(q_1),T^{-n}(q_2),\ldots)
q_0^nq_1^{n^2/2}q_2^{n^3/3}\cdots .
$$
\end{thm}
\begin{proof} First note
\begin{equation}T(q_0)^nT(q_1)^{n^2/2}T(q_2)^{n^3/3}\cdots =
q_0q_1^{-\half}q_2^{\frac{1}{3}}\cdots q_0^{n-1}q_1^{(n-1)^2/2}\cdots.
\end{equation}
Write $\Omega = \sum_{-\infty}^\infty V_n(q_1,\ldots)q_0^nq_1^{n^2/2}\cdots$. Then
\begin{multline*}q_0q_1^{-\half}q_2^{\frac{1}{3}}\cdots\sum_nV_{n-1}(q)q_0^{n-1}q_1^{(n-1)^2/2}\cdots
= \\
q_0q_1^{-\half}q_2^{\frac{1}{3}}\cdots\Omega(q) = \Omega(T(q)) = \\
\sum V_n(T(q))T(q_0)^nT(q_1)^{n^2/2}\cdots = \\
q_0q_1^{-\half}q_2^{\frac{1}{3}}\cdots\sum_nV_{n}(T(q))q_0^{n-1}q_1^{(n-1)^2/2}\cdots.
\end{multline*}
Taking coefficients of $q_0^{n-1}$ we get $V_n(T(q))=V_{n-1}(q)$. The formula follows since 
$V_0=V$. 
\end{proof}
\begin{exam} Consider the formula in the theorem with $q_j\mapsto 1,\ j\ge 2$. We have
$$T(q)|_{q_j\mapsto 1,j\ge 2} = q_1
$$
and by \eqref{012}
$$V(q_1,1,1,\ldots) = q_1^{-\frac{1}{24}}\sum_\lambda q_1^{p_1(\lambda)} =
q_1^{-\frac{1}{24}}\prod_{m\ge 1}(1-q_1^m)^{-1}=\eta(q_1)^{-1}.
$$
The assertion of the theorem is then
\begin{equation}\eta(q_1)^{-1}\sum_{n\in \Z} q_0^nq_1^{n^2/2} = q_1^{\frac{-1}{24}}\prod_{m\ge
0}(1+q_0q_1^{m+\half})(1+q_0^{-1}q_1^{m+\half}).
\end{equation}
Multiplying through by $\eta(q_1)$ yields the triple product formula
\eqref{22}.
\end{exam}

\section{Quasimodular forms}\label{sec:quasi}

In this section, we recall the theory of quasimodular forms as developed in \cite{KZ}. 
All results are due to Kaneko and Zagier and are recalled here solely for the convenience of the
reader. 

We fix a subgroup of finite index $\Gamma\subset\Gamma_1 := \text{\rm SL}(2,\Z)$. A holomorphic modular
form of weight $k$ (for $\Gamma$) is a holomorphic function $f(\tau)$ on the upper half-plane
$\sH=\{\tau=x+iy\ |\ y>0\}$ satisfying
\begin{equation}\label{31}f\left(\frac{a\tau+b}{c\tau+d}\right) = (c\tau+d)^kf(\tau),\quad \begin{pmatrix}a & b
\\ c & d\end{pmatrix}\in\Gamma\,.
\end{equation}
Note for some $\ell\ge 1$
we have $(\begin{smallmatrix}1 & \ell \\ 0 & 1\end{smallmatrix})\in\Gamma$, so $f$ may be expanded
in a Fourier series
\begin{equation}\label{32}f(\tau) = \sum a_n\exp(2\pi in\tau/\ell) = \sum a_n q^{n/\ell}\,.
\end{equation}
We assume that $f$ is holomorphic at $i\infty$, i.e. $a_n=0$ for $n<0$ and the Fourier series
converges for $|q|<1$. The holomorphic modular forms constitute a graded ring
$$M_*(\Gamma) := \oplus M_k(\Gamma)
$$
graded by the weight $k$. 

An {\it almost holomorphic} function on $\sH$ will be a function
\begin{equation} F(\tau) = \sum_{m=0}^N f_m(\tau)y^{-m}
\end{equation}
on $\sH$ where $y=\text{\rm Im}(\tau)$ and each $f_m$ has a Fourier expansion as in \eqref{32}. An
almost holomorphic modular form of weight $k$ is an almost holomorphic function $f(\tau)$
satisfying the weight $k$ modularity property \eqref{31}. 
\begin{exam} The classical Eisenstein series ($B_k$ Bernoulli number as in Lemma (\ref{lem22}))
\begin{equation}\label{34}G_k(q):=\frac{-B_k}{2k}+\sum_{n=1}^\infty\Big(\sum_{d|n}
d^{k-1}\Big)q^n,\quad k=2,4,6,\ldots
\end{equation}
is modular of weight $k$ for $\Gamma=\Gamma_1$ and $k\ge 4$. On the other hand, $G_2$ satisfies
the transformation
\begin{equation}\label{35} G_2\left(\frac{a\tau+b}{c\tau+d}\right) =
(c\tau+d)^2G_2(\tau)-\frac{c(c\tau+d)}{4\pi i} \,.
\end{equation}
Notice, however, that
\begin{equation}\label{36} \left[\text{\rm Im}\left(\frac{a\tau+b}{c\tau+d}\right)\right]^{-1} -(c\tau+d)^2y^{-1} =
-2ic(c\tau+d) \,. 
\end{equation}
It follows that the function
\begin{equation}\label{37} G_2^*(\tau) := G_2(\tau)+\frac{y^{-1}}{8\pi}
\end{equation}
is an almost holomorphic modular form of weight $2$ for $\Gamma_1$. 
\end{exam}
\begin{defn}\label{defn32} A quasimodular form of weight $k$ is a holomorphic function $f(\tau)$ on
$\sH$ such that there exists an almost holomorphic modular form $F = \sum_{m=0}^N f_my^{-m}$ of
weight $k$ with $f_0=f$. 
\end{defn}
\begin{exam} $G_2$ is quasimodular of weight $2$. Indeed, one can take $F=G_2^*$. 
\end{exam}

Let us write $AHM_*(\Gamma)$ (resp. $QM_*(\Gamma)$) for the graded ring of almost holomorphic
(resp. quasimodular) forms. (Kaneko and Zagier write $\widehat{M}$ and $\widetilde{M}$, but this
makes it difficult to remember which is which.)
\begin{prop}\label{prop34} The assignment $F=\sum_{j=0}^My^{-j}f_j\mapsto f_0$ defines an
isomorphism of graded rings, $AHM_*(\Gamma)\cong QM_*(\Gamma)$.
\end{prop}
\begin{proof} Note this is well defined, i.e. $\sum_{j=0}^My^{-j}f_j(\tau) \equiv 0$ for
holomorphic $f_j$ if and only if all the $f_j$ are zero. To see this one can e.g.\ apply the
differential operator $\frac{iM}{2}y^{-1}+\frac{d}{d\bar\tau}$ and argue by induction on $M$. The
map in question is surjective by definition, so it suffices to show injectivity. Suppose for some
$r\ge 1$ and $f_r \ne 0$ that
$F=y^{-r}f_r+\ldots+y^{-M}f_M$ is almost holomorphic modular of weight $k$. Let
$A=(\begin{smallmatrix}a & b \\ c & d\end{smallmatrix}) \in\Gamma$ and write $j=c\tau+d$. Using
\eqref{36} we get
\begin{multline*} (j^2y^{-1}-2icj)^rf_r(A\tau)+\ldots +(j^2y^{-1}-2icj)^Mf_M(A\tau) = \\
j^ky^{-r}f_r(\tau)+\ldots+j^ky^{-M}f_M(\tau).
\end{multline*}
Now identify coefficients of powers of $y^{-1}$ 
\begin{align}\label{39} f_M(A\tau) &= j^{k-2M}f_M(\tau) \,,\\
f_{M-1}(A\tau)-2icj\binom{M}{1}f_M(A\tau) &= j^{k-2M+2}f_{M-1}(\tau) \notag \,,\\
f_{M-2}(A\tau)-2icj\binom{M-1}{1}f_{M-1}(A\tau) & +(2icj)^2\binom{M}{2}f_M(A\tau) \notag\\
&= j^{k-2M+4}f_{M-2}(\tau)\notag\,, \\
\makebox[5cm][c]{$\vdots$} & \makebox[3cm][c]{$\vdots$} \notag 
\displaybreak[0]\\
f_r(A\tau)-2icj\binom{r+1}{1}f_{r+1}(A\tau) \qquad\quad& \notag \\ 
+(2icj)^2\binom{r+2}{2}f_{r+2}(A\tau) -\ldots &=
j^{k-2r}f_r(\tau)\notag \,.
\end{align}
The terms not involving $y$ give
\begin{equation}(-2icj)^rf_r(A\tau)+\ldots+(-2cij)^Mf_M(A\tau) = 0. \label{310}
\end{equation} 
Solve \eqref{39} recursively:
\begin{align*}f_M(A\tau) &= j^{k-2M}f_M(\tau) \,,\\
f_{M-1}(A\tau) &= 2ic\binom{M}{1}j^{k-2M+1}f_M(\tau)+j^{k-2M+2}f_{M-1}(\tau)\notag \,,\\
f_{M-2}(A\tau) &= \ldots + j^{k-2M+4}f_{M-2}(\tau) \notag \,, \\
\makebox[2cm][c]{$\vdots$} & \makebox[5cm][c]{$\vdots$} \notag \\
f_r(A\tau) & = \ldots + j^{k-r}f_r(\tau)\,.  \notag
\end{align*}
Finally, substituting in \eqref{310} yields 
\begin{multline}\label{312} (c\tau+d)^{k-r}c^rf_r(\tau)
+\alpha_{r+1}(c\tau+d)^{k-r-1}c^{r+1}f_r(\tau) +\ldots \\
+ \alpha_M (c\tau+d)^{k-M}c^Mf_M(\tau) = 0.
\end{multline}
Here the $\alpha_j$ are constants independent of $c$ and $d$. Varying $A\in \Gamma$ yields a
contradiction. Indeed, the map $\Gamma\to\C^2,\ A\mapsto (c,d)$ has Zariski dense image, so in the
above identity, $c$ and $d$ can be taken to be independent variables. The coefficient of
$c^rd^{k-r}$ is nontrivial. 
\end{proof}
\begin{prop}\begin{enumerate}\item[(i)] $M_*(\Gamma)\subset QM_*(\Gamma)$.
\item[(ii)] $QM_*(\Gamma) = M_*(\Gamma)\otimes\C[G_2]$.
\item[(iii)] $QM_*(\Gamma)$ is stable under the operator (of degree $2$) $D:= \frac{d}{d\tau}$. We
have 
$$QM_k(\Gamma) = \oplus_{0\le i\le k/2}D^iM_{k-2i}(\Gamma)\oplus\C\cdot D^{k/2-1}G_2.
$$
\end{enumerate}
\end{prop}
\begin{proof} (i) is clear. To prove (ii), we claim first that the map
\begin{equation}\label{313}M_*(\Gamma)\otimes\C[G_2^*] \to AHM_*(\Gamma)
\end{equation}
is an isomorphism, where $G_2^*$ is as in \eqref{37}. Indeed, for $F = f_0+\ldots+y^{-M}f_M$ almost
holomorphic modular of weight $k$, it follows from the first line of \eqref{39} that $f_M$ is
holomorphic modular of weight $k-2M$. We have
$$F - f_M\cdot (8\pi G_2^*)^M = g_0+\ldots+y^{-M+1}g_{M-1}
$$
is almost holomorphic modular of weight $k$. Surjectivity of \eqref{313} follows by induction on
$M$. Injectivity is straightforward, keeping track of powers of $y^{-1}$. Assertion (ii) now
follows from proposition (\ref{prop34}). Finally, (iii) is left for the reader.
\end{proof}

\section{Quasimodularity for characters $\Omega$ and $V$}\label{sec:qmov}

The purpose of this section is to prove
\begin{thm} The series $\Omega(\tau_0,\tau_1,\ldots)$ and $V(\tau_1,\tau_2,\ldots)$ are 
quasimodular of weights $0$ and $-\half$ respectively. 
\end{thm}
We focus first on $\Omega$. Recall (definition (\ref{defn32})) a series
$F(\tau_0,\tau_1,\tau_2,\ldots)$ is said to be quasimodular of weight
$k$ if it can be expanded in a formal Taylor series
$$\sum_{J=(j_0,j_2,j_3,\ldots)} B_J(\tau_1)\tau^J/J!
$$
with $B_J(\tau_1)$ quasimodular of weight 
\begin{equation}\label{4wt} k+\wt(J):= k+j_0+3j_2+4j_3+\ldots. 
\end{equation}
Since we have not
specified a group $\Gamma$, there will be no harm e.g. in replacing $\tau_1$ by $2\tau_1$. 
Define
$$\Phi(q)=\Omega(-q_0,q_1,q_2,\ldots).$$ 
One has $\Omega(q) =
\Phi(q_0^2,q_1^2,q_2^2,\ldots)/\Phi(q)$, so, since the space of quasimodular forms is a ring, it
will suffice to prove
$\Phi(q)$ is quasimodular. We have
$$\Phi(1,x_1,1,1,\ldots) = x_1^{-\frac{1}{24}}\prod_{r\ge
0}(1-x_1^{r+\frac{1}{2}})^2 = \big(\eta(x^{1/2})/\eta(x)\big)^2,
$$
which is modular of weight $0$. It therefore suffices to show
$$G :=\log\big(\Phi(x)/\Phi(1,x_1,1,1,\ldots)\big)
$$
is quasimodular. We have (with $\xi(s)$ as in \eqref{16})
\begin{multline}\label{41}
G = -2\pi i\big(\sum_{n\ge 1}\xi(-2n-1)\tau_{2n+1}\big)+ \\
\sum_{r\ge
0}\big[\log(1-x_0x_1^{n+\frac{1}{2}}x_2^{(n+\frac{1}{2})^2}\cdots)+
\log(1-x_0^{-1}x_1^{n+\frac{1}{2}}x_2^{-(n+\frac{1}{2})^2}\cdots)- \\
2\log(1-x_1^{n+\frac{1}{2}})\big]\,.
\end{multline}
We expand the final sum in \eqref{41} 
\begin{multline*}B:= \\ -\sum_{n\ge 0\,,\ l\ge 1}\frac{x_1^{(n+\half)l}}{l} 
\left(\exp\left(2\pi
il\left(\tau_0+\left(n+\thalf\right)^2\tau_2+\left(n+\thalf\right)^3\tau_3+\ldots\right)\right)\right. \\
+\left.
\exp\left(2\pi
il\left(-\tau_0-\left(n+\thalf\right)^2\tau_2+\left(n+\thalf\right)^3\tau_3-\ldots\right)\right) - 2\right).
\end{multline*}
Define coefficients $\alpha_m(J)$, where $J=(j_0,j_2,j_3,\ldots)$ and $\ j_k \ge 0$, by
\begin{equation}\frac{\partial^J}{\partial \tau^J}B|_{\tau_0=\tau_2=\ldots=0} =
\sum_{m=1}^\infty\alpha_m(J)x_1^{m/2}\,.
\end{equation}
Note $\alpha_m(0)=0$. For $J\ne 0$, write $|J|=j_0+j_2+\ldots$. We get
\begin{equation}\alpha_m(J) = -(2\pi i)^{|J|}\sum_{n\ge 0\,,\ l\ge
1}l^{|J|-1}(1+(-1)^{\wt(J)})\left(n+\thalf\right)^{\wt(J)-|J|}x_1^{(n+\half)l}.
\end{equation}
Thus, $\alpha_m(J)=0$ for $\wt(J)$ odd. For even weight, we find
\begin{multline}\label{45}
\frac{\partial^J}{\partial \tau^J}B|_{\tau_0=\tau_2=\ldots=0} = -2(2\pi
i)^{|J|}2^{|J|-\wt(J)}\sum_{m\ge 1}m^{|J|-1}x_1^{m/2}\sum_{d|m;\ d \text{\rm  odd}}d^{\wt(J)-2|J|+1}
\\ = -2(2\pi i)^{|J|}2^{|J|-\wt(J)}\sum_{m\ge 1}m^{\wt(J)-|J|}x_1^{m/2}\sum_{l|m;\ \frac{m}{l}
\text{\rm  odd}}l^{2|J|-\wt(J)-1} \,. 
\end{multline}

The Eisenstein series of weight $k$ and level $1$ was given in \eqref{34}. 
As in \cite{KZ}, we work with the level two Eisenstein series
\begin{gather*}F_k^{(1)}(q) := G_k(q^{\half}) - G_k(q) =
\sum_{n=1}^\infty\Big(\sum_{d|n,2\not\,\mid d}(n/d)^{k-1}\Big)q^{n/2}\,, \\
F_k^{(2)} := G_k(q^{\half}) - 2^{k-1}G_k(q) = (1-2^{k-1})\zeta(1-k)/2+
\sum_{n=1}^\infty\Big(\sum_{d|n,2\not\,\mid d}d^{k-1}\Big)q^{n/2}.
\end{gather*}
Suppose first that the inequality $k(J) := \wt(J)-2|J|+2>1$ holds. It follows from the first
equality in \eqref{45} that
\begin{equation}\label{48}
\frac{\partial^J}{\partial \tau^J}B|_{\tau_0=\tau_2=\ldots=0} =
-4\pi
i(2^{-k(J)+1})(\partial/\partial\tau_1)^{|J|-1}(F_{k(J)}^{(2)}(x_1)-F_{k(J)}^{(2)}(0))\,.
\end{equation}
Note $k(J)$ is even, so the remaining possibility is $k(J) \le 0$. In this case, the second
equality in \eqref{45} implies 
\begin{equation}\frac{\partial^J}{\partial \tau^J}B|_{\tau_0=\tau_2=\ldots=0} = -2(2\pi
i)^{2-k(J)}(\partial/\partial\tau_1)^{\wt(J)-|J|}(F_{2-k(J)}^{(1)}(x_1))\,.
\end{equation}
Since the $F^{(i)}_k$ are quasimodular, it follows that the $\frac{\partial^J}{\partial
\tau^J}B|_{\tau_0=\tau_2=\ldots=0}$ are quasimodular except possibly (because of the constant
term on the right) in in the case $k(J)\ge 2;\ |J|=1$ where the constant term is not correct. 
This only happens if $J=(0,\ldots,0,1,0,\ldots)$ where $j=1$ for some $j\ge 3$ odd and all
the other entries are zero. In this case $k(J)=j+1$ and it follows from \eqref{48} and
the identity $\xi(s)=(2^s-1)\zeta(s)$ that
\begin{multline*}\frac{\partial}{\partial\tau_j} B|_{\tau=0} = (-4\pi
i)2^{-j}\big(F^{(2)}_{j+1}(x_1)+(2^j-1)\zeta(-j)/2\big)= \\
(-4\pi i)2^{-j}F^{(2)}_{j+1}(x_1)+2\pi i\xi(-j)\,.
\end{multline*}
The last constant term exactly cancels the $-2\pi i\xi(-j)$ which appears in the first sum on
the right in \eqref{41}. This completes the proof of quasimodularity for $\Omega$. 

It remains to prove quasi-modularity for $V(x_1,\ldots)$. The proof parallels that
in \cite{KZ}. We begin by expanding
\begin{gather} V(\tau_1,\tau_2,\ldots) = \sum_{K =
(k_2,k_3,\ldots)}A_K(\tau_1)\tau^K/K!\,, \\
V(T^{-n}(\tau)) = \sum_{m=0}^\infty\sum_{K}(\partial/\partial\tau_1)^mA_K(\tau_1)
(T^{-n}(\tau_1)-\tau_1)^m(m!)^{-1}(T^{-n}(\tau))^K/K!\,, \\
\label{413}
\Omega(\tau) =
\sum_{n=-\infty}^\infty\sum_{J=(j_0,j_2,j_3,\ldots)}\sum_{m=0}^\infty\sum_{K} 
(\partial/\partial\tau_1)^mA_K(\tau_1)(T^{-n}(\tau_1)-\tau_1)^m(m!)^{-1}\times \\
\times (T^{-n}(\tau))^K(K!)^{-1}a_Jn^{\text{\rm wt}(J)}\tau^J(J!)^{-1}x_1^{n^2/2}\,. \notag
\end{gather} 
Here $\sum_Ja_Jn^{\text{\rm wt}(J)}\tau^J(J!)^{-1}$ is the power series expansion of
$x_0^nx_2^{n^3/3}x_3^{n^4/4}\cdots$, with $a_J$ independent of $n$. 

One verifies
the expansion
\begin{equation}\label{414}
T^{-n}(\tau_j) =
\tau_j+\binom{j+1}{1}n\tau_{j+1}+\binom{j+2}{2}n^2\tau_{j+2}+\ldots\quad
.
\end{equation}
Fix $m\ge 0$ and $K=(k_2,k_3,\ldots)$. Let $P(\tau)=\tau_2^{p_2}\tau_3^{p_3}\ldots$
and $Q_i(\tau)=\tau_i^{q_{ii}}\tau_{i+1}^{q_{i,i+1}}\ldots;\ i\ge 2$ be monomials
with
$\deg(P)=\sum p_j = m$ and $\deg(Q_i)=k_i$. Define the weights by $\wt(P) = \sum
(j+1)p_j$ and $\wt(Q_i) = \sum_{j\ge i} (j+1)q_{ij}$. As a consequence of
\eqref{414}, $P$ appears in $(T^{-n}(\tau_1)-\tau_1)^m(m!)^{-1}$ with a
coefficient $c_{P,m}n^{\wt(P)-2m}$, where $c_{P,m}$ is constant independent
of $n$. Similarly, $Q:=\prod Q_i$ appears in $(T^{-n}(\tau))^K(K!)^{-1}$ with
coefficient $c_{\{Q_i\},K}n^{\sum_i \wt(Q_i)-(i+1)k_i}$. Now fix
$J=(j_0,j_2,\ldots)$. We obtain a contribution to the coefficient of the monomial
$P(\tau)Q(\tau)\tau^J$ in \eqref{413} of the form
\begin{equation}\label{415}
c_{J,K,m,\{Q_i\},P}\frac{\partial^m}{\partial\tau_1^m}A_K(\tau_1)
\sum_{n=-\infty}^\infty n^{\wt(P)-2m+\sum_i (\wt(Q_i)-(i+1)k_i)+\wt(J)}x_1^{n^2/2}.
\end{equation}
If the exponent of $n$ is odd, this cancels. Assume this exponent equals $2r$
for $r\ge 0$ an integer. Then \eqref{415} can be rewritten
\begin{equation}\label{416}
b_{J,K,m,\{Q_i\},P}\frac{\partial^m}{\partial\tau_1^m}A_K(\tau_1)
\frac{\partial^r}{\partial\tau_1^r}\theta_{00}(\tau_1)\,,
\end{equation}
where $\theta_{00}(\tau_1)=\sum \exp(\pi in^2\tau_1)$ is a modular form of weight
$1/2$. Since differentiation preserves quasi-modularity and increases weight by
$2$, it follows that, assuming $A_K(\tau_1)$ is quasimodular of weight
$\wt(K)-\half$, the expression in \eqref{416} is quasi-modular of weight  
$\wt(P)+\sum\wt(Q_i)+\wt(J)$.

We will prove $A_K(\tau_1)$ is quasimodular of weight $\wt(K)-1/2$ by induction on
$\wt(K)$. Write
\begin{multline*}\theta_{00}(x) = \sum_{-\infty}^\infty x^{n^2/2}= \prod_{m\ge
1}(1-x^m)\prod_{m\ge 0}(1+x^{m+\half})^2 = \\
\eta(x)\Omega(1,x,1,1\ldots)\,.
\end{multline*}
Substituting $x_j=1;\ j\ne 1$ in Theorem (\ref{thm23}) yields
\begin{gather*}\Omega(1,x_1,1,\ldots) = V(x_1,1,\ldots)\sum x_1^{n^2/2} =
V(x_1,1,\ldots)\theta_{00}(x_1) \,,\notag \\
V(x_1,1,\ldots) = \eta(x_1)^{-1} \,. 
\end{gather*}
This is modular of weight $-1/2$ as desired. Now fix an index set
$M=(m_2,m_3,\ldots)$ with $\wt(M)>0$ and assume
$A_K(\tau_1)$ is quasimodular of weight $\wt(K)-1/2$ for all $K$ with
$\wt(K)<\wt(M)$. Consider the coefficient $B_M(\tau_1)$ of $\tau^M$ in the
expansion of $\Omega(\tau)$ \eqref{413}. We know that $B_M$ is quasimodular of
weight $\wt(M)$. It is a sum of terms \eqref{416}, of which all but one involve
$A_K$ with $\wt(K)<\wt(M)$ and are quasimodular of weight $\wt(M)$ by our inductive
hypothesis. The one remaining term, which appears with coefficient $(M!)^{-1}$, is
$A_M(\tau_1)\theta_{00}(\tau_1)$ (take $J=0=m,\ P=1,\ Q_i = \tau_i^{m_i}$.)
It follows that $A_M(\tau_1)$ is quasimodular of weight $\wt(M)-1/2$ as desired.

\section{Preliminaries on partitions}\label{sec:prepart}
Let $\lambda = \lambda_1\ge \lambda_2\ge\ldots\ge\lambda_\ell >0=\lambda_{\ell+1}=\ldots$ be a
partition. $\ell=\ell(\lambda)$ is the length of the partition, and $|\lambda| = \sum\lambda_i$
is the number being partitioned. 

Let $f(\lambda)$ be a function on the set of partitions. Assuming $f$ does not grow too
rapidly with $|\lambda|$ we may consider for $|q|\ll 1$ the ratio 
\begin{equation}\label{51} \langle f\rangle_q\ := \frac{\sum_\lambda f(\lambda)q^{|\lambda|}}{\sum_\lambda
q^{|\lambda|}} = (q)_\infty \sum_\lambda f(\lambda)q^{|\lambda|}.
\end{equation}
Here, we write for $q$ given and $n\le \infty$, 
\begin{equation}\label{52}(a)_n = (1-a)(1-aq)(1-aq^2)\cdots (1-aq^n).
\end{equation}
Fix an integer $n\ge 1$ and variables $t_1,\ldots,t_n$. We will be interested in the following
``$n$-point correlation function''
\begin{equation}\label{53} F(t_1,\ldots,t_n) := \Big<\prod_{k=1}^n\Big(\sum_{i=1}^\infty
t_k^{\lambda_i-i+\half}\Big)\Big>_q.
\end{equation} 

To understand the relation between $F$ and the character $V(q)$, we take a small detour, beginning
with some basic ideas from the theory of partitions. The {\it diagram} of a partition 
$\lambda$ is simply the array of dots (or squares) with
$\lambda_1$ dots in the first row, $\lambda_2$ dots in the second row, etc. Here is $4,4,3,2,1$.
$$\begin{array}{cccc} \bullet & \bullet & \bullet & \bullet \\
\bullet & \bullet & \bullet & \bullet \\
\bullet & \bullet & \bullet \\
\bullet & \bullet  \\
\bullet
\end{array}
$$
The transposed partition $\lambda'$ is obtained by flipping along the diagonal
$$\begin{array}{ccccc} \bullet & \bullet & \bullet & \bullet & \bullet\\
\bullet & \bullet & \bullet & \bullet \\
\bullet & \bullet & \bullet \\
\bullet & \bullet 
\end{array}
$$
For $\lambda_i\ge i$ let $m_i=\lambda_i-i$ be the number of dots to the right of the $i$-th point
on the diagonal, and let $n_i$ be the number of points below it. (so
$n_i(\lambda)=m_i(\lambda')$.) Note $m_1>m_2>\ldots \ge 0$ and similarly for the $n_i$. The $m_i$
and $n_i$ are called Frobenius coordinates of $\lambda$ and written
$(m_i,\ldots,m_p|n_1,\ldots,n_p)$. For example, Frobenius coordinates for $4,4,3,2,1$ are
$(3,2,0|4,2,0)$. Frobenius coordinates are related to the $\lambda_i$
via the following generating
function identity.
\begin{lem}\label{lem51} Let $\lambda$ be a partition of length $\ell$ with Frobenius
coordinates\newline
$(m_i,\ldots,m_p|n_1,\ldots,n_p)$. Then
$$\sum_i (t^{\lambda_i-i+\half}-t^{-i+\half}) = \sum_k (t^{m_k+\half}-t^{-(n_k+\half)}).
$$
\end{lem}
\begin{proof} It is clear that the $m_k$ correspond to the nonnegative $\lambda_i-i$. Thus it will
suffice to identify the terms with negative exponents on the two sides. This follows from the
identity of sets of numbers
$$\{\l_{p+1}-(p+1),\dots,\l_{\ell}-\ell,-n_1-1,\dots,-n_p-1\} = \{-\ell,\dots,-2,-1\},,
$$ 
which, in turn, follows by comparing the negative elements from  the identity
$$
\{\l_1-1,\dots,\l_{\ell}-\ell,-\l'_1,\dots,-\l'_k +
k-1\}=\{-\ell,\dots, k-1\}\,, \forall k \ge \l_1\,,
$$
established in \cite{MacD}, Chapter I, (1.7).
\end{proof}

Substituting $t=\exp(u)$ and expanding in the previous lemma yields identities
\begin{multline}\label{54} \sum_i \left(\lambda_i-i+\thalf\right)^r
  +(-1)^{r+1} \left(i-\thalf\right)^r =
 \\
\sum_j \left(m_j+\thalf\right)^r + (-1)^{r+1}\left(n_j+\thalf\right)^r =: p_r(\lambda).
\end{multline}
For example, $|\lambda| = p_1(\lambda)$. Recall we have already encountered the $p_r(\lambda)$ in
the expansion for the character $V(q)$:
\begin{equation} V(q_1,q_2,q_3,\ldots) = q_1^{-\xi(-1)}q_3^{-\xi(-3)}\cdots \sum_\lambda
\prod_{r\ge 1}q_r^{p_r(\lambda)}.
\end{equation}
To see the connection more specifically, write $q_n = \exp(2\pi i \tau_n)$ and consider the
differential operator 
\begin{equation}\label{56} \delta = \delta(u) := u^{-1}+(2\pi i)^{-1}\sum_{n=1}^\infty
\frac{u^n}{n!}\frac{\partial}{\partial \tau_n}.
\end{equation}

\begin{prop} We have with $F$ as in \eqref{53} and $t=\exp(u)$
$$F(t) = \eta(q_1)\delta V|_{\tau_2=\tau_3=\ldots =0}.
$$
\end{prop}
\begin{proof} Consider the Bernoulli polynomials, defined by
\begin{equation}\label{57} \frac{u\exp(xu)}{\exp(u)-1} = \sum_{k=0}^\infty B_k(x)u^k/k!;\quad
B_0(x) = 1,\ B_1(x) = x-\half
\end{equation}
One has (\cite{hida}, theorem 1, p. 43)
\begin{equation}\label{58} \xi(1-m) = \frac{-B_m(\half)}{m}.
\end{equation}
Combining these identities yields
\begin{equation}\label{59} -\sum_{n=1}^\infty
\xi(-n)u^n/n! = 
\frac{1}{t^\half - t^{-\half}}-u^{-1} = \sum_{i=1}^\infty t^{-i+\half} -u^{-1}.
\end{equation}
On the other hand, using lemma \ref{lem51} and the identity
\begin{equation}\label{510} \eta(q_1)^{-1} = q_1^{-\xi(-1)}\sum_\lambda q_1^{|\lambda|},
\end{equation}
one computes
\begin{multline*}\eta(q_1)\delta V|_{\tau_2=\ldots = 0} = u^{-1} -\sum \xi(-n)u^n/n! + \\
+ \Big(\sum_\lambda q_1^{|\lambda|}\Big)^{-1}\sum_\lambda q_1^{|\lambda|}\sum_{i=1}^\infty
(t^{\lambda_i - i+\half} - t^{-i+\half} ) = \\
= u^{-1}+\Big(\sum_\lambda q_1^{|\lambda|}\Big)^{-1}\sum_\lambda q_1^{|\lambda|}\sum_{i=1}^\infty
t^{\lambda_i - i+\half} - u^{-1} = \\
=\Big<\sum_{i=1}^\infty t^{\lambda_i - i+\half}\Big>_{q_1} = F(t).
\end{multline*}
\end{proof}
More generally, one finds
\begin{equation}\label{512} \eta(q_1)\delta(u_1)\circ\ldots\circ\delta(u_n)V|_{\tau_2=\ldots=0} 
= \Big<\prod_{k=1}^n\Big(\sum_{i=1}^\infty
t_k^{\lambda_i-i+\half}\Big)\Big>_{q_1} = F(t_1,\dotsc,t_n).
\end{equation}
This is also equivalent to the formula \eqref{016a} from the introduction.

\section{The formula for correlation functions}\label{sec:results}
In this section we give a detailed statement of our results on the generating function
$F(t_1\dotsc,t_n)$ defined in \eqref{53}. To simplify notation, we
write  $\T(x)=\T(x;q)$ for the
following theta function
\begin{align}\label{61}
\T(x) = \To(x;q)&=\eta^{-3}(q)\sum_{n\in{\Bbb Z}} (-1)^n q^{\frac{(n+1/2)^2}2} x^{n+1/2} \\
&=\qi^{-2} (x^{1/2}-x^{-1/2}) (q x\i (q/x\i \,. 
\end{align}
The function $\T(x)$ is odd,
\begin{equation*}
\T(x^{-1})=-\T(x), 
\end{equation*}
and satisfies the difference equation
\begin{equation}\label{62}
\T(q^m x)=(-1)^m q^{-m^2/2} x^{-m} \T(x)\,,
\quad m\in\Z \,. 
\end{equation}
Define
\begin{equation}
\T^{(k)}(x):= \left(x\d x\right)^k \To(x;q)\, ,\ \T' = \T^{(1)},\ \text{\rm etc.}
\end{equation}
where $x\d x$ is the natural invariant vector
field on the group $\C^* = \C\setminus 0$.
We shall prove the following

\begin{thm}\label{thm61}
\begin{equation}\label{64}
\vspace{-3 \jot}
F(\tdn)= 
\sum_{\sigma\in\mathfrak{S}(n)}\, \frac 
{\displaystyle \det\left( \frac{\displaystyle \T^{(j-i+1)}(t_{\sigma(1)}\cdots
    t_{\sigma(n-j)})}{\displaystyle (j-i+1)!} 
\right)_{i,j=1}^n}
{\displaystyle \T(t_{\sigma(1)})\, \T(t_{\sigma(1)} t_{\sigma(2)}) \dots
\T(t_{\sigma(1)}\cdots t_{\sigma(n)})}
\end{equation} 
Here $\sigma$ runs through all permutations $\mathfrak{S}(n)$ of $\{1,\dotsc,n\}$,the
matrices in the numerator have size $n\times n$, and we define 
$1/(-n)!=0$ if $n\ge 1$. 
\end{thm} 
In particular, since $\T'(1)=1$ and $\T''(1)=0$, we have 
\begin{align*}
F(t_1) &=\frac{1}{\T(t_1)}\,,\notag\\
F(t_1,t_2) &=\frac{1}{\T(t_1 t_2)}
\left(
\frac{\T'(t_1)}{\T(t_1)} + \frac{\T'(t_2)}{\T(t_2)} \right)\,. \notag 
\end{align*}
For $n=3$ the formula \eqref{64} looks as follows:
\begin{multline}\label{66}
F(t_1,t_2,t_3)=\\
\frac1{\displaystyle\T(t_1 t_2 t_3)}
\sum_{\sigma\in\mathfrak{S}(3)}
\det
\left(
 \begin{array}{rrc}
  \frac{\displaystyle\T'\left(t_{\sigma(1)} t_{\sigma(2)}\right)}
      {\displaystyle \T\phantom{{}'}\left(t_{\sigma(1)} t_{\sigma(2)}\right)}&
  \frac{\displaystyle 1}{\displaystyle 2} \,
  \frac{\displaystyle\T''\left(t_{\sigma(1)}\right)}
    {\displaystyle\T\phantom{{}''}\left(t_{\sigma(1)}\right)}&
  \frac{\displaystyle\T'''(1)}{\displaystyle 3!}\\
  {\displaystyle 1}&
  \frac{\displaystyle\T'\phantom{{}'}\left(t_{\sigma(1)}\right)}
  {\displaystyle\T\phantom{{}''}\left(t_{\sigma(1)}\right)}& 
  0\\
  &
  {\displaystyle 1}&
  \T'(1)
 \end{array}
\right)\\
= \frac{1}{\T(t_1 t_2 t_3)}
\left( \sum_{1\le i \ne j \le 3}
\frac{\T'(t_i)}{\T(t_i)} 
\frac{\T'(t_i t_j)}{\T(t_i t_j)} -
\sum_{i=1}^3 \frac{\T''(t_i)}{\T(t_i)} + \T'''(1)
\right) \,.
\end{multline}
The determinants in \eqref{64}, having only one non-zero diagonal
below the main diagonal, 
have a nice combinatorial expansion, see
the formula \eqref{78} below. 

One has the following simple
\begin{lem}\label{lem61a}
\begin{equation}\label{th'}
\T^{(2m+1)}(1;q)=(2m+1)! 
\sum_{k_1+2k_2+\dots=m} \frac{(-2)^{k_1+k_2+\dots}}
{k_1!\, k_2!\, \cdots} 
\prod_i \left(\frac{G_{2i}(q)}{(2i)!}\right)^{k_i} \,.
\end{equation}
\end{lem}
Observe that the sum in \eqref{th'} is over all 
partitions $1^{k_1} 2^{k_2} 3^{k_3} \dots$ of 
the number $m$ and $k_1+k_2+\dots$ is the length
of such a partition. In particular,  
\begin{align*}
\T^{(3)}(1;q)&= - 6 G_2(q)\,,\\ 
\T^{(5)}(1;q)&= - 10 G_4(q)+ 60 G_2(q)^2\,, \\
\T^{(7)}(1;q)&= - 14 G_6(q)+420 G_4(q) G_2(q)- 840 G_2(q)^3\,. 
\end{align*}
This lemma will be proved below.

\begin{remark} The equality in Theorem \ref{thm61} may look like an 
equality of multivalued functions but in fact
the only ambiguity is the factor
$$
\sqrt{\,t_1\cdots t_n}
$$
which appears on both sides of of the theorem and,
therefore, can be ignored.   
\end{remark}
\begin{remark}  Alternatively, one may consider the $(n+1)\times (n+1)$ matrix
\begin{equation}\label{67}\left(\begin{array}{ccccc} \T(t_1\cdots t_n) & \T'(t_1\cdots t_{n-1}) &\hdots &
\frac{1}{n!}\T^{(n)}(1) \\
0 & \T(t_1\cdots t_{n-1}) & \hdots & \frac{1}{(n-1)!}\T^{(n-1)}(1) \\
\vdots & \vdots & \hdots & \vdots \\
0 & 0 & \hdots & \T'(1) \\
0 & 0 & \hdots & 0
\end{array}\right)
\end{equation}
(Note the bottom right hand entry is $\T(1)=0$.) This matrix has rank $n$, so there will be a
unique column vector 
$$\begin{pmatrix}v_n(\td{n}) \\
v_{n-1}(\td{n-1}) \\
\vdots \\
v_1(t_1) \\
1\end{pmatrix}
$$
which is killed by multiplication on the left by the matrix \eqref{67}. Then
$$F(\td{r}) = \sum_{\sigma} v_r(t_{\sigma(1)},\dotsc,t_{\sigma(r)}).
$$
\end{remark}

\begin{proof}[Proof of Lemma \ref{lem61a}]
Let $f(a)$ be some function of a variable $a$ and let $b$ be 
another variable. We have:
\begin{align*}
e^{f(a+b)-f(a)}&=\exp\left(\sum_{m=1}^\infty \frac{f^{(m)}(a)}{m!}
  b^m\right)\\
&=\prod_{m=1}^\infty \exp\left(\frac{f^{(m)}(a)}{m!}
  b^m\right)\\
&=\prod_{m=1}^\infty \sum_{k=0}^\infty \left(\frac{f^{(m)}(a)}{m!}
  \right)^k \frac{b^{mk}}{k!}\\
&=\sum_{s=0}^\infty
\sum_{k_1+2k_2+3k_3+\dots=s}
\frac{b^s}{k_1!\, k_2!\, k_3!\cdots}
\prod_{i=1}^\infty\left(\frac{f^{(i)}(a)}{i!}\right)^{k_i}\,.
\end{align*}
In other words, we have 
\begin{multline}\label{67a}
e^{-f(a)} \frac{d^s}{d a^s} e^{f(a)} = \\
s! \sum_{k_1+2k_2+3k_3+\dots=s} \frac1{k_1!\, k_2!\, k_3!\cdots}
\left(\frac{f'}{1!}\right)^{k_1} 
\left(\frac{f''}{2!}\right)^{k_2}
\left(\frac{f'''}{3!}\right)^{k_3} \cdots \,.
\end{multline} 
We wish to apply the above formula to the 
triple product formula for
the theta function and use the fact that
$$
\left.\left(x\frac{d}{dx}\right)^{m}
\log\left((qx)_\infty (q/x)_\infty\right)
\right|_{x=1}= 
\begin{cases} 
-2 G_{m}(q) - B_{m}/m\,, &\text{$m$ is even}\\
0\,, &\text{$m$ is odd} \,.
\end{cases} 
$$
The derivative $\T^{(2m+1)}(1;q)$ is quasi-modular of weight $2m$.
It is, therefore, sufficient, to compute the term
$$
\left. (2m+1) \left(x\frac{d}{dx} (x^{1/2}-x^{-1/2}) \right)
\left(x\frac{d}{dx}\right)^{2m} (qx)_\infty (q/x)_\infty 
\right|_{x=1}\,,
$$
which gives the contribution of the maximal weight. Using 
\eqref{67a} we obtain that the above term equals
the RHS of \eqref{th'} modulo terms of lower weight.
This establishes \eqref{th'}.
\end{proof}

Before diving into the proof of Theorem \ref{thm61}
in the next section, we give as a warmup the proof for the $1$-point
function.
\begin{thm}\label{thm64}
$$\Big<\sum_i t^{\l_i -i+\half}\Big>_q = \frac{1}{\T(t)}.
$$
\end{thm}
\begin{lem}\label{lem65}
$$\Big<t^{\l_k}\Big>_q = \frac{(q)_{\infty}}{(1-q)\cdots (1-q^{k-1})(1-q^kt)(1-q^{k+1}t)\cdots }
$$
\end{lem}
\begin{proof} The notation is as in \eqref{52}. Recall in section (\ref{sec:prepart}) we defined
the transpose $\l'$ of a partition $\l$. The lemma is straightforward once one observes that
$$\Big<t^{\l_k}\Big>_q = \Big<t^{\l_k'}\Big>_q = \Big<t^{\#\{i|\l_i\ge k\}}\Big>_q.
$$
\end{proof}
Recall the following Heine's $q$-analog of the Gauss
${}_2F_1$-summation. Given any $a$, $b$, $c$ satisfying 
$|c|<|ab|$ we have 
\begin{equation}\label{69}
\sum_{n=0}^\infty \frac{(a)_n (b)_n}{(c)_n (q)_n} 
\left(\frac{c}{ab}\right)^n = \frac{(c/a)_\infty
    (c/b)_\infty}{(c)_\infty (c/ab)_\infty}\,,
\end{equation}
see for example \cite{GR}, Section 1.5. Recall that the symbol 
$(a)_n$ was defined in \eqref{52}. 

\begin{proof}[Proof of theorem \ref{thm64}] 
\begin{multline*} \Big<\sum_i t^{\l_i -i+\half}\Big>_q = \sum_{k=1}^\infty t^{\half
-k}\frac{(q^k)_\infty}{(q^kt)_\infty} = \\
= t^{-\half}\frac{(q)_\infty}{(qt)_\infty}\sum_{r=0}^\infty
t^{-r}\frac{(qt)_i}{(q)_i} = \frac{(q)_\infty^2}{(t^\half -
t^{-\half})(qt)_\infty(q/t)_\infty}= \frac{1}{\T(t)}.
\end{multline*}
\end{proof}

\section{Beginning of the proof of Theorem \ref{thm61}}\label{sec:not}
Our strategy for proving theorem \ref{thm61} will be based on the
following simple implication
$$
\Bigg(
\begin{matrix}
\text{\rm $f(x)$ is holomorphic on $\C\setminus 0$ and }\\
f(qx)=f(x)
\end{matrix}
\Bigg)
\Rightarrow
\Bigg(
\text{\rm $f(x)$ is a constant}
\Bigg)\quad . 
$$
We shall show that both sides of theorem \ref{thm61} satisfy
the same difference equation and have the 
same singularities. 

We now introduce some useful notation.   
Given indices 
$
1< \iln \,.
$
set
\begin{align}\label{71}
\ex{i_1,\dots,i_r}{t_1,\dots,t_r}
:&= \la \prod_{k=1}^n  t_k^{\l_{i_k}-i_k+1/2}
\ra \\ 
&= \prod_{k=1}^n  t_k^{1/2-i_k}
\la \prod_{k=1}^n  t_k^{\l'_{i_k}}
\ra \notag \\
&= \qi
\frac
{t_1^{1/2-i_1} \dots t_r^{1/2-i_r}}
{(q)_{i_1-1} (q^{i_1} t_1)_{i_2-i_1} (q^{i_2} t_1 t_2)_{i_3-i_2}
\dots
(q^{i_r} t_1 \cdots t_r)_\infty}
\,,\notag
\end{align}

The notation $(a)_n$ is as in \eqref{52}, and the last equality in \eqref{71} is a
straightforward generalization of lemma (\ref{lem65}). For brevity, write
$$
\eit:=\ex{i_1,\dots,i_r}{t_1,\dots,t_r} \,.
$$
Set 
\begin{align}
H(\tdn)&:=\sum_{\iln} \eit \,, \label{72} \\
G(\tdn)&:=\sper H(\tsdn) \,. \label{73}
\end{align}

By $\P n$ denote the set of all partitions
of the 
set $\on$. An element $\pi$ of $\P n$ 
$$
\pi=\{\pi_1,\dots,\pi_\ell\} \in \P n\,, \quad \pi_i\subset\on\,,
$$
is by definition an unordered collection of subsets 
of $\on$ such that
\begin{align*}
&\pi_i \cap \pi_j = \emptyset\,, \quad i\ne j\,, \\
&\pi_1 \cup \dots \cup \pi_\ell = \on \,.
\end{align*}
Note the difference between partitions of the
\emph{set} $\on$ and partitions of the \emph{number} $n$. 
The number $\ell$ is called the \emph{length} of $\pi$ and is denoted
by $\lpi$. The subsets $\pi_i$ are called the \emph{blocks} of $\pi$.

Given a partition $\pi\in\P n$ set
$$
G^\pi(t_1,\dots,t_n):=G\left(\prod_{k\in\pi_1} t_k\,, \dots , 
\prod_{k\in\pi_{\lpi}} t_k \right) \,.
$$
For example, if $n=3$ and $\pi=\{\{1,2\},\{3\}\}$ then
$$
G^\pi(t_1,t_2,t_3)=G(t_1 t_2, t_3) \,.
$$
With this notation we have
\begin{equation}\label{76}
F(\tdn)=\sum_{\pi\in\P n} G^\pi(t) \,. 
\end{equation}
Set also
\begin{equation}\label{77}
F^\pi(\tdn):= F\left(\prod_{k\in\pi_1} t_k\,, \dots , 
\prod_{k\in\pi_{\lpi}} t_k \right)\,, \quad \pi\in\P n \,.
\end{equation}

We shall also need the set $\G n$ of all {\it 
compositions} of the set $\on$. By definition,
the set $\G n$ consists of all {\it ordered}
collections
$$
\g=(\g_1,\dots,\g_\ell)\,, \quad \g_i\subset\on\,,
$$
such that
\begin{align*}
&\g_i \cap \g_j = \emptyset\,, \quad i\ne j\,, \\
&\g_1 \cup \dots \cup \g_\ell = \on \,.
\end{align*}
The number $\ell$ is called again the {\it length} of $\g$ and is denoted
by $\lga$.  The set $\G n$ is in a natural one-to-one
correspondence with the set $\U_n$ of all partial flags 
$$
\emptyset=\ups_1\varsubsetneq\ups_2\varsubsetneq\dots\varsubsetneq\ups_\ell
\varsubsetneq\ups_{\ell+1}=\on
$$
in the set $\on$, namely, we set
\begin{equation}\label{77a}
\ups(\g)_k=\g_1 \cup \dots \cup \g_{k-1}\,.
\end{equation} 
For any partition $\pi$ or composition $\g$ of an $n$-element  set of
length $\ell$ define
its sign as 
$$
(-1)^\pi = (-1)^\gamma = (-1)^{n+\ell}\,.
$$
Then, in particular, the sign of a permutation equals the sign of
the corresponding partition into disjoint cycles.

Denote by 
$
U(\tdn)
$
the RHS of \eqref{64}, and set
\begin{multline*}
T(\tdn):= \T(t_1 \dotsc t_n) U(\tdn) \\
=\sum_{\sigma\in\mathfrak{S}(n)}
\frac 1{\T(t_{\sigma(1)})\T(t_{\sigma(1)} t_{\sigma(2)}) \dotsc \T(t_{\sigma(1)}\cdots
t_{\sigma(n-1)})}\times \notag \\
\times \det\left( \frac{\T^{(j-i+1)}(t_{\sigma(1)}\cdots t_{\sigma(n-j)})}{(j-i+1)!} 
\right)\,. \notag
\end{multline*}
Expanding the determinants one obtains
\begin{equation}\label{78}
T(\tdn)=\sum_{\g\in\G n} (-1)^{\g}
\T^{(\#\g_1)}(1) \,
\prod_{k=2}^{\lga}
\frac{\T^{(\#\g_k)}\left(\prod_{i\in\ups(\g)_{k}} t_i\right)}
{\T\left(\prod_{i\in\ups(\g)_{k}} t_i\right)}\,, 
\end{equation}
where we used the notation \eqref{77a} for the partial sums of a 
composition $\g$. 
Define the functions $T^\pi(t)$, $\pi\in\P n$, as in \eqref{77}. 

We will use two subsets of $\P n$. By 
$\Po n$ denote the set of all partitions $\pi$ such
that $\pi$ has at most one  block of cardinality $>1$ which,
in addition, contains the number 1. 
By $\Pt n$ denote the set of all partitions $\pi$ 
that have $\{1\}$ as a block. We have, for example,
\begin{gather*}
\{\{1,2\},\{3\}\}\in\Po 3\,,\\
\{\{1\},\{2,3\}\}\in\Pt 3\,.
\end{gather*}
Denote by 
$$
\text{Atom}({\on})=\{\{1\},\dotsc,\{n\}\} \in \P n
$$
the partition into 1-element subsets. It is clear
that
$$
\Po n \cap \Pt n = \text{Atom}({\on})\,.
$$
Define the subsets $\Go n,\Gt n \subset \G n$ in the same way as for
partitions of $\on$.

We shall need the following identities 

\begin{lem}
\begin{gather}\label{710}
\sum_{\pi\in\P n} (-1)^{\pi} \lpi !  =
\sum_{\g\in\G n} (-1)^{\g} = 1  \\
\label{711} \sum_{\pi\in\P n} (-1)^{\lpi} (\lpi-1) !  = 0 \,.
\end{gather}
\end{lem}
\begin{proof} We shall prove \eqref{710}. The proof of \eqref{711} is similar. Define $f: \prod_n
\to \prod_{n-1}$ in the obvious way, by simply omitting $n$ from a partition. We argue by
induction on $n$. Note 
$$\# f^{-1}(\pi) = \ell(\pi)+1.
$$
More precisely, $f^{-1}(\pi)$ contains $\lpi$ partitions of length $\lpi$ and one
partition of length $\lpi+1$. We compute
\begin{multline*} \sum_{\pi\in\prod_n} (-1)^{n+\lpi}\lpi! = \\
\sum_{\pi\in\prod_{n-1}}
(-1)^n\Big((-1)^{\lpi} \lpi \lpi!+(-1)^{\lpi+1}(\lpi+1)!\Big) = \\
\sum_{\pi\in\prod_{n-1}}
(-1)^{n-1+\lpi}\lpi! = 1.
\end{multline*}
\end{proof}

\section{Difference equations for the correlation functions}\label{sec:difF}
The aim of this section is the following

\begin{thm}\label{thm81}
$$
F(q t_1, t_2,\dots,t_n)= -q^{1/2} t_1\cdots t_n \left(
\sum_{\pi\in\Po n} (-1)^{\pi} F^\pi (\tdn) \right) \,.
$$
\end{thm}

Recall that, by definition,
\begin{multline*}
\sum_{\pi\in\Po n} (-1)^{\pi} F^\pi (\tdn)=\\
\sum_{s=0}^{n-1} 
(-1)^s 
\sum_{1<i_1<\dots < i_s \le n} F(t_1 t_{i_1} t_{i_2} \cdots t_{i_s},
\dots,\widehat{\,t_{i_1}}, \dots,\widehat{\,t_{i_s}}, \dots) \,.\notag
\end{multline*}
We begin with some auxiliary propositions. 

\begin{prop}\label{prop82} For any $k=1,\dots,n$

\begin{align*}
H(t_1,\dots,qt_k,\dots,t_n)=
&-q^{1/2} t_1 \cdots t_n  H(\tdn) \notag\\
&+ \frac{q t_k}{1-q t_k} H(t_1,\dots,t_{k-1},
qt_k t_{k+1}, t_{k+2},\dots,t_n) \\
&- \frac{1}{1-q t_k} H(t_1,\dots,t_{k-2},
qt_{k-1} t_{k}, t_{k+1},\dots,t_n) \,.\notag
\end{align*}

Here if $k=1$ (or $k=n$) the third (second) summand
should be omitted. 
\end{prop}

\begin{cor}\label{cor83}
\begin{align}
G(q t_1, t_2,\dots,t_n)=& -q^{1/2} t_1\cdots t_n  G(\tdn) \label{82} \\
&-\sum_{k=2}^n G(q t_1 t_k,t_2,\dots,t_{k-1},t_{k+1},\dots, t_n)\notag \\
=& -q^{1/2} t_1\cdots t_n  \left(
\sum_{\pi\in\Po n} (-1)^{\pi} (n-\lpi)! \, \, G^\pi (t) \right) \,. \label{83}
\end{align}
\end{cor}

\begin{lem}\label{lem84} Let $u$ and $v$ be two variables.
For any integers $a < b$ 
\begin{equation}
\sum_{a < i < b} \frac{(q v)^{1/2-i}}
{(q^a u)_{i-a} (q^i u v)_{b-i}} = 
\frac{1}{1-q v}
\left(
\frac{(q v)^{3/2-b}}{(q^a u)_{b-a-1}} -
\frac{(q v)^{1/2-a}}{(q^{a+1} u v )_{b-a-1}}
\right) 
\,.
\end{equation}
In particular, for $a=1$ and $u=1$ we obtain
\begin{equation}
\sum_{1 \le  i < b} \frac{(q v)^{1/2-i}}
{(q)_{i-1} (q^i v)_{b-i}}
=
\frac{1}{1-q v}
\frac{(q v)^{3/2-b}}{(q)_{b-2}} \,. 
\end{equation}
\end{lem}

\begin{proof} Follows from the following particular
case of the $q$-Gauss summation formula \eqref{69}: for any $c$ and $d$ one
has 
\begin{equation}
\sum_{i=c}^{\infty} (q v)^{1/2 - i}
\frac{(q^d u v)_{i-c}}{(q^d u)_{i-c}}=
 -(q v)^{3/2-c} \frac{1-q^{d-1} u}{1- q v} \,.
\end{equation}
\end{proof}
In the rest of the section we shall write  $a\ll b$ if $a+1<b$.
\begin{lem}\label{lem85} We have
\begin{align}
&H(t_1,\dots,qt_k,\dots,t_n)= \notag \\
&\quad=
\sum_{i_1<\dots < i_{k} \ll i_{k+1} < \dots < i_n}
q^{1/2-i_k} (1-q^{i_k} t_1 \cdots t_{k}) t_{k+1} \cdots t_n
\ex{i_1,\dots,i_n}{t_1,\dots,t_n}\label{88} \\
&\quad=
\sum_{i_1<\dots < i_{k-1} \ll i_{k} < \dots < i_n}
q^{3/2-i_k} (1-q^{i_k-1} t_1 \cdots t_{k-1}) t_k \cdots t_n
\ex{i_1,\dots,i_n}{t_1,\dots,t_n} \label{89}
\end{align}
\end{lem}
\begin{proof} 
Recall that by definition
$$
H(t_1,\dots,qt_k,\dots,t_n) = \sum_{\iln} 
\ex{i_1,\dots,i_n}{t_1,\dots,qt_k,\dots,t_n} \,.
$$
We have the two following identities: 
\begin{align}
&\ex{i_1,\dots,i_n}{t_1,\dots,qt_k,\dots,t_n} = \notag \\
&\qquad=
q^{1/2-i_k} (1-q^{i_k} t_1 \cdots t_{k}) t_{k+1} \cdots t_n
\ex{i_1,\dots,i_k,i_{k+1}+1,\dots,i_n+1}{t_1,\dots,t_n} \label{810} \\
&\qquad=
q^{1/2-i_k} (1-q^{i_k} t_1 \cdots t_{k-1}) t_k \cdots t_n
\ex{i_1,\dots,i_{k-1},i_{k}+1,\dots,i_n+1}{t_1,\dots,t_n} \label{811}
\end{align}

Replacing in \eqref{810} the indices  $i_{k+1}+1,\dots,i_n+1$ 
by $i_{k+1},\dots,i_n$ we obtain \eqref{88}. 
Replacing in \eqref{811} the indices  $i_{k}+1,\dots,i_n+1$ 
by $i_{k},\dots,i_n$ we obtain \eqref{89}. 
\end{proof}

\begin{proof}[Proof of proposition \eqref{prop82}]
Consider the general case $1 < k < n$. The two other cases $k=1$ and $k=n$ are
similar. Rewrite the sum \eqref{88} as follows
\begin{multline}\label{812}
\sum_{i_1<\dots < i_{k} \ll i_{k+1} < \dots < i_n}
q^{1/2-i_k} (1-q^{i_k} t_1 \cdots t_{k}) t_{k+1} \cdots t_n
\ex{i_1,\dots,i_n}{t_1,\dots,t_n} \\
 = 
\sum_{i_1< \dots < i_n} (\dots)  - 
\sum_{\stackrel{i_1<\dots < i_{k} \ll i_{k+2} < \dots < i_n}{
i_{k+1}=i_k+1}}  (\dots)\,, 
\end{multline} 
and denote the two summands in the RHS of \eqref{812} by $\s1$ and $\s2$. The
formula \eqref{88} together with the following identity
\begin{multline*}
q^{1/2-i_k} (1-q^{i_k} t_1 \cdots t_{k}) t_{k+1} \cdots t_n
\ex{i_1,\dots,i_k,i_k+1,i_{k+2},\dots,i_n}{t_1,\dots,t_n} = \\
q^{1/2-i_k} (1-q^{i_k} t_1 \cdots t_{k+1}) t_{k+2}  \cdots t_n
\ex{i_1,\dots,i_{k-1},i_k,&i_{k+2},\dots,i_n}{t_1,\dots,t_{k-1},t_k t_{k+1},&t_{k+2},
\dots,t_n} 
\end{multline*}
implies that
$$
\s2=H(t_1,\dots,q t_k t_{k+1}, \dots, t_n) \,.
$$
Set
\begin{align}
\s3:&=\s1 + q^{1/2} t_1\cdots t_n H(\tdn) \notag \\
&=t_{k+1}\cdots t_n \sum_{\iln} q^{1/2-i_k} {\eit}\ \,. \label{814}
\end{align}
Using Lemma \ref{lem84} we can evaluate one of the
nested sums in \eqref{814} as follows:
\begin{multline*}
\sum_{i_{k-1} < i_k < i_{k+1}} 
\frac {(qt_k)^{1/2-i_k}} 
{(q^{i_{k-1}} t_1 \cdots t_{k-1})_{i_k-i_{k-1}}
(q^{i_{k}} t_1 \cdots t_{k})_{i_{k+1}-i_{k}}} \\= 
\frac1{1-q t_k} 
\left(
\frac{(qt_k)^{3/2-i_{k+1}}}
{(q^{i_{k-1}} t_1 \cdots t_{k-1})_{i_{k+1}-i_{k-1}-1}} -
\frac{(qt_k)^{1/2-i_{k-1}}}
{(q^{i_{k-1}+1} t_1 \cdots t_{k})_{i_{k+1}-i_{k-1}-1}}
\right) \,.
\end{multline*}
Then sum $\s3$ becomes the following difference
\begin{align*}
&\frac{t_k \cdots t_n}{1-q t_k} 
\sum_{\dots<i_{k-1} \ll i_{k+1} < \dots}  
q^{3/2-i_{k+1}} (1- q^{i_{k+1}-1}t_1\cdots t_{k-1})
\ex{\dots,i_{k-1},i_{k+1},\dots}{
\dots,t_{k-1}, t_k t_{k+1},\dots}
\\
&-\frac{t_{k+1} \cdots t_n}{1-q t_k}
\sum_{\dots<i_{k-1} \ll i_{k+1} < \dots}  
q^{1/2-i_{k-1}} (1- q^{i_{k-1}}t_1\cdots t_{k})
\ex{\dots,i_{k-1},i_{k+1},\dots}{
\dots,t_{k-1}t_k, t_{k+1}, \dots} \,.
\end{align*} 
By Lemma \ref{lem85} this yields
\begin{multline*}
\s3= \frac1{1-qt_k} H(t_1,\dots,t_{k-1},q t_k t_{k+1}, \dots, t_n)
\\
- \frac1{1-qt_k} H(t_1,\dots,q t_{k-1}t_k, t_{k+1}, \dots, t_n) \,.
\end{multline*}
Thus,
\begin{align*}
H(t_1,\dots,qt_k,\dots,t_n)=& \s3 - 
q^{1/2} t_1 \cdots t_n H(\tdn) -\s2 \\
=&-q^{1/2} t_1 \cdots t_n  H(\tdn)\\
&+ \frac{q t_k}{1-q t_k} H(t_1,\dots,t_{k-1},
qt_k t_{k+1}, t_{k+2},\dots,t_n) \\
&- \frac{1}{1-q t_k} H(t_1,\dots,t_{k-2},
qt_{k-1} t_{k}, t_{k+1},\dots,t_n) \,. 
\end{align*}
This concludes the proof. 
\end{proof}

\begin{proof}[Proof of corollary \ref{cor83}]
We will prove \eqref{82}; equation \eqref{83} follows easily from \eqref{82}.
Substitute the formula from proposition \ref{prop82} into \eqref{73} and consider the
coefficient of a summand
$$
H(\dots, q t_1 t_j, \dots)\,, \quad j=2,\dots,n \,.
$$
Such a summand arises from the expansion of
$$
H(\dots, q t_1, t_j, \dots)\quad\text{\rm\ and}
\quad H(\dots, t_j, q t_1,  \dots) \,.
$$
The sum of the coefficients equals
$$
\frac{q t_1}{1- q t_1} -
\frac{1}{1- q t_1} = -1 \,.
$$
This proves \eqref{82}. 
\end{proof}

\begin{prop}\label{prop86}
\begin{equation}\label{815}
F(q t_1, t_2,\dots,t_n)= -q^{1/2} t_1\cdots t_n \left(
\sum_{\pi\in\Pt n} G^\pi(t) \right) \,. 
\end{equation}
\end{prop}
\begin{proof}
Substitute \eqref{83} into \eqref{76}. Let $c_\pi$ denote the
coefficient of the summand $G^\pi(t)$ in this 
expansion. It is clear that
$$
c_\pi = -q^{1/2} t_1\cdots t_n \,, \pi \in \Pt n \,.
$$
Suppose that $\pi \notin \Pt n$ and show that $c_\pi=0$.
Since $F(q t_1, t_2,\dots,t_n)$ 
is symmetric in variables $t_2,\dots,t_n$
we can assume that
$$
\pi=\{\{1,2,\dots,m\},\pi_2,\dots,\pi_l\}\,, \quad 1<m\le n\,,
$$
where $\{\pi_2,\dots,\pi_l\}$ is a partition of the set
$\{m+1,\dots,n\}$.  

Let $\si$ be a partition of $\{1,\dots,m\}$
$$
\si=\{\si_1,\dots,\si_s\} \in \P m \,.
$$
Set
$$
\pi\land\sigma=\{\si_1,\dots,\si_s,\pi_2,\dots,\pi_l\} \in \P n \,.
$$
It is easy to see that the term  $G^\pi(t)$ arises precisely
from the expansion of the terms of the form
$$
G^{\pi\land\sigma} (qt_1,t_2,\dots,t_n)\,, \quad \si\in\P m \,.
$$
By \eqref{83} it arises with the coefficient
$$
q^{1/2} t_1\cdots t_n (-1)^{\ell(\si)} (\ell(\si)-1)! \,.
$$
By \eqref{711}, the sum of these coefficients over all
$\si\in\P m$ equals zero. 
This concludes the proof.
\end{proof}

\begin{proof}[Proof of theorem \ref{thm81}]
We will show that the RHS in theorem \ref{thm81} equals the RHS of \eqref{815}. Substitute
\eqref{76} into RHS of theorem \ref{thm81} and let $c_\pi$ denote the
coefficient of the summand $G^\pi(t)$ in this 
expansion. Again, it is clear that
$$
c_\pi = -q^{1/2} t_1\cdots t_n \,, \pi \in \Pt n \,.
$$
and we have to show that
$$
c_\pi = 0 \,, \pi \notin \Pt n \,.
$$
Again, we can suppose that
$$
\pi=\{\{1,2,\dots,m\},\pi_2,\dots,\pi_l\}\,, \quad 1<m\le n\,,
$$
where $\{\pi_2,\dots,\pi_l\}$ is a partition of the set
$\{m+1,\dots,n\}$.

For any set $S$ denote by $\text{Atom}(S)$ the partition of $S$ into
one-element blocks. The summand $G^\pi(t)$ arises in the RHS of
from the expansion of terms
$$
F^{\left\{\{S\},\text{Atom}(\on\setminus S)\right\}}(t)\,, 
\quad  1\in S \subset \{1,\dots,m\}\,,
$$
and therefore the coefficient $c_\pi$ equals
\begin{align*}
c_\pi&=-q^{1/2} t_1\cdots t_n \sum_{1\in S \subset \{1,\dots,m\}} (-1)^{|S|+1} \\
&=-q^{1/2} t_1\cdots t_n \sum_{k=0}^{m-1} (-1)^k \binom{m-1}k \\
&=0 \,.
\end{align*} 
This concludes the proof. 
\end{proof}

\section{Singularities of correlation functions}\label{sec:sing}
In this section we consider the singularities of the 
function $F(\tdn)$. The series \eqref{76} converges 
uniformly on compact subsets of the set
$$
\left\{|t_i|>1\,, i=1,\dots,n\right\}\setminus
\bigcup_{S\subset\on} \bigcup_{m=1}^\infty
\left\{q^m \prod_{k\in S} t_k =1 \right\}
$$
and has simple poles on the divisors
$$
\left\{q^m \prod_{k\in S} t_k =1 \right\}\,,
\quad S\subset\on \,.
$$
The difference equation from theorem \ref{thm81} gives a meromorphic
continuation of the function $F(t)$ onto 
(a double cover of) the domain $({\Bbb C}\setminus 0)^n$. 

By symmetry, it suffices to consider the divisors
\begin{equation}\label{91}
q^m t_1 t_2 \dots t_k = 1\,, \quad k\le n,\quad  m=1,2,\dots \,. 
\end{equation}
We shall prove the following 

\begin{thm}\label{thm91} We have 
\begin{equation}\label{92}
F(t_1,\dots,t_n)=(-1)^{m}
 \frac{q^{m^2/2} m^{k-1}}{q^m t_1\dots t_k-1}
\frac{F(t_{k+1},\dots,t_n)}{(t_{k+1}\dots t_n)^m}+ \dots \,, 
\end{equation}
where dots stand for terms regular at the divisor \eqref{91} and we
assume that 
$$
F(t_{k+1},\dots,t_n)=1\,,\quad k=n\,.
$$
\end{thm}

A curious property of \eqref{92} is that the residue 
\begin{equation}\label{93}
\Res_{\{q^m t_1 t_2 \dots t_k = 1\}} F(\tdn) 
\end{equation}
does not depend on the variables
$t_1,\dots,t_k$.

Introduce the following notation. Given a function
$f(i_1,\dots,i_n)$ set
$$
\tsm f(i_1,\dots,i_n) = 
\sum_{1\le i_1 \le \dots \le i_n \le m}
\frac{f(i_1,\dots,i_n)}
{\# \Stab_{S(n)} \{i_1,\dots,i_n\}} \,.
$$
If the function $f$ is symmetric, we have
\begin{equation}\label{94}
\tsm f(i_1,\dots,i_n) = \frac 1{n!} \sum_{i_1,\dots,i_n=1}^m
f(i_1,\dots,i_n) \,. 
\end{equation}
With this notation we have
$$
F(t)=
\sum_{\stackrel{\text{\rm over all permutations}}{
\text{\rm of $t_1,\dots,t_n$}}}
 \,\,
\widetilde{\sum_{1\le i_1 \le \dots \le i_n}}
\eit \,.
$$
It is clear that the evaluation of the residue \eqref{93} boils down to 
the evaluation of the following sum
\begin{multline}\label{95}
\sum_{\stackrel{\text{\rm permutations}}{\text{\rm of $t_1,\dots,t_k$}}} \\
\tsm \frac{t_1^{1/2-i_1} \dots t_k^{1/2-i_k}}
{(q)_{i_1-1} (q^{i_1} t_1)_{i_2-i_1} (q^{i_2} t_1 t_2)_{i_3-i_2} \dots
(q^{i_k} t_1 \dots t_k)_{m-i_k}} \,, 
\end{multline}
where the variables $t_1,\dots,t_k$ are subject to 
constraint
$$
q^m t_1 \cdots t_k =1 \,.
$$
This evaluation follows from a curious rational function identity
which shall now be established. Set
$$
\cbm {i_1,\dots,i_k}{t_1,\dots,t_k}:=
\frac1
{(q)_{i_1-1} (q^{i_1} t_1)_{i_2-i_1} (q^{i_2} t_1 t_2)_{i_3-i_2} \dots
(q^{i_k} t_1 \dots t_k)_{m-i_k}} \,.
$$
Then we have the following

\begin{thm}\label{thm92} Let $m,k\in{\Bbb N}$ and suppose that
$$
q^m t_1 t_2 \dots t_k = 1\,.
$$
Then
\begin{equation}\label{96}
\scy\quad \tsm \cbm {i_1,\dots,i_k}{t_1,\dots,t_k} =
\frac{m^{k-1}}{(k-1)!} \,. 
\end{equation}
In particular, this sum is independent of variables $q,t_1,\dots,t_k$.
\end{thm}

\begin{proof}[Proof of Theorem \ref{thm92}]
Let us use $q,t_1,\dots,t_{k-1}$ as independent variables. 
First, show that LHS of \eqref{96} does not depend on $t_1,\dots,t_{k-1}$. Since
it is a rational function bounded on infinity it suffices
to show that it does not have any poles. The only poles it can 
possibly have are simple poles on divisors
$$
t_s t_{s+1} \dots t_{r-1} t_r = q^{-l}\,, \quad s\le r\le k-1,\,
1\le l \le m \,.
$$
By symmetry, it suffices to consider the pole on
\begin{equation}\label{97}
t_1 \dots t_r = q^{-l}\,, \quad \quad r\le k-1,\,
1\le l \le m \,. 
\end{equation}
This pole arises from two type of  summands in \eqref{96}. First, it
arises from the summands
\begin{equation}\label{98}
\cbm {i_1,\dots,i_k}{t_1,\dots,t_k}\,, \quad i_r \le l < i_{r+1}\,, 
\end{equation}
and it also arises from summands
\begin{equation}\label{99}
\cbm {i_1,\dots,i_k}{t_{r+1},\dots,t_k,t_1,\dots,t_r}\,, 
\quad i_{k-r} \le m-l < i_{k-r+1}\,.
\end{equation}
We shall match each summand of the form \eqref{98} to a
summand of the form \eqref{99} in such a way that the poles
in each pair will cancel out. Namely, it is easy to see
that the sum 
\begin{equation}\label{910}
\cbm {i_1,\dots,i_k}{t_1,\dots,t_k} \,+\, 
\cbm {i_{r+1}-l,\dots,i_k-l,i_1+m-l,\dots,i_r+m-l}
{t_{r+1},\dots,t_k,t_1,\dots,t_r}\,, 
\end{equation}
where
$$
1\le i_1\le \dots \le i_r \le l < i_{r+1}\le \dots \le i_k \le m\,,
$$
is regular on the divisor \eqref{97}; notice that the
inequalities
$$
1\le i_1 \le \dots \le i_r < i_{r+1} \le \dots \le i_k \le m
$$
imply that
$$
\#\Stab\{i_1,\dots,i_k\}=
\#\Stab\{i_{r+1}-l,\dots,i_k-l,i_1+m-l,\dots,i_r+m-l\}\,.
$$
This proves that the LHS of \eqref{96} is regular on \eqref{97}.

Therefore, the LHS of \eqref{96} does not depend on $t_1,\dots,t_k$
and we can substitute
$$
t_1=t_2=\dots=t_{k-1}=1\,,\quad t_k=q^{-m} \,.
$$
Then the LHS of \eqref{96} becomes 
\begin{equation}\label{911}
\tsm \,\,
\sum_{r=1}^k  \frac1{(q)_{i_r} (q^{i_r-m})_{m-i_r}} \,. 
\end{equation}
Since the inner sum in \eqref{911} is symmetric in $i_1,\dots,i_k$
we can use \eqref{94} to obtain
\begin{align*}
&\frac1{k!} \sum_{i_1,\dots,i_k=1}^m \sum_{r=1}^k 
\frac1{(q)_{i_r} (q^{i_r-m})_{m-i_r}} \notag\\
&= \frac{m^{k-1}}{(k-1)!} \sum_{i=1}^m 
\frac1{(q)_{i} (q^{i-m})_{m-i}} \notag\\
&= \frac{m^{k-1}}{(k-1)!} \frac1{(q^{1-m})_{m-1}}
\sum_{i=1}^m 
\frac{(q^{1-m})_{i-1}}{(q)_{i}} \notag\\
&=\frac{m^{k-1}}{(k-1)!} \,,\notag
\end{align*}
where the last equality follows from the q-binomial theorem.
This concludes the proof. 
\end{proof}

Replacing in \eqref{96} the variables $q,t_1,\dots,t_k$ by their
reciprocals one obtains the following

\begin{cor}\label{cor93} Suppose that $q^m t_1 t_2 \dots t_k = 1$. Then
the sum
\begin{equation}\label{912}
\scy
\tsm \frac{t_1^{1/2-i_1} \dots t_k^{1/2-i_k}}
{(q)_{i_1-1} (q^{i_1} t_1)_{i_2-i_1} (q^{i_2} t_1 t_2)_{i_3-i_2} \dots
(q^{i_k} t_1 \dots t_k)_{m-i_k}} 
\end{equation}
equals
$$
(-1)^{m-1} \frac{q^{m^2/2} m^{k-1}}{(k-1)!}\,,
$$
for all $m$ and $k$.
\end{cor}

This corollary gives the evaluation of the sum \eqref{95} and
this proves \eqref{92}

\section{Difference equations for the RHS in \eqref{64}}\label{sec:dif}

In this section we show that the functions $U(\tdn)$ satisfy the 
very same difference equations as the functions  $F(t)$ do.

Since we shall use only the difference equation \eqref{62} for
the theta function $\T(x)$, let us consider following 
general situation. 

Let $t_0$ be an auxiliary variable. We shall eventually
let
$$
t_0\to 1 \,.
$$
Modify temporarily the definition \eqref{77a} as
follows
\begin{equation}\label{101}
\ups(\g)_k=\{0\}\cup \g_1 \cup \dots \cup \g_{k-1}\,.
\end{equation}
 
Suppose that a  function
$$
r(x;m)\,, \quad m=0,1,2,\dots\,,
$$
satisfies the two following properties:
\begin{equation}\label{102}
r(x;0)=1 
\end{equation}
identically and
\begin{equation}\label{103}
r(qx;m)=\sum_{i=0}^m (-1)^i \binom mi\, r(x;m-i) \,. 
\end{equation}
Consider the function
\begin{equation}\label{104}
R(\tdn\ve t_0)=\sum_{\g\in\G n} (-1)^{\g} \prod_{k=1}^{\lga}
r\left({\tprod}_{i\in\ups(\g)_{k}} t_i \,; \#\g_k
\right) \,,
\end{equation}
where we use the modified definition \eqref{101}. We want to show that
\begin{thm}\label{thm101} Suppose a function $r(x;m)$ satisfies \eqref{102}
and \eqref{103}. Then 
\begin{equation}\label{105}
R(qt_1,t_2,\dots,t_n\ve t_0)= \sum_{\pi\in\Po n} 
(-1)^{\pi} R^\pi(\tdn\ve t_0) \,, 
\end{equation}
where the function $R(\tdn\ve t_0)$ was defined in \eqref{104}
\end{thm}

Here, as usual, 
\begin{equation}\label{106}
R^\pi(\tdn\ve t_0):= R\left(\left.\prod_{k\in\pi_1} t_k\,, \dots , 
\prod_{k\in\pi_{\lpi}} t_k \,\right|\, t_0 \right)\,, \quad \pi\in\P n \,.
\end{equation}

\begin{proof}  Given a composition 
$$
\g=(\g_1,\dots,\g_l)\,,
$$
let $m=m(\g)$ denote the number of the block that contains 1.
Substitute \eqref{103} into \eqref{104}.
We obtain
\begin{equation}\label{107}
R(qt_1,t_2,\dots,t_n\ve t_0)=
\sum_{\g\in\G n} 
\sum_{s_{m+1}=0}^{\#\g_{m+1}}
\dots
\sum_{s_l=0}^{\#\g_{l}} 
\sqb 
{\g_{m+1},\dots,\g_l}
{s_{m+1},\dots,s_l}\,, 
\end{equation}
where 
\begin{equation}\label{108}
\sqb
{\g_{m+1},\dots,\g_l}
{s_{m+1},\dots,s_l}
\end{equation}
stands for the following product
\begin{multline*}
(-1)^{n+\lga+s_{m+1}+\dots+s_{l}}
\prod_{k=1}^{m}
r \left({\tprod}_{i\in\ups(\g)_{k}} t_i\,; \#\g_k
\right)\times \\
\prod_{k=m+1}^{l} \binom {\#\g_i}{s_i} 
r \left({\tprod}_{i\in\ups(\g)_{k}} t_i\,; \#\g_k - s_i 
\right)
\,. 
\end{multline*}

Let us divide all summands \eqref{108} into 3 following types
according to 
occurrence of certain patters in the sequence $s_{m+1},\dots,s_l$.
The summands of the form
\begin{equation}\label{109}
\sqb
{\g_{m+1},&\dots,&\g_k,&\dots}
{\#\g_{m+1},&\dots,&\#\g_k,&0,\dots,0}\,, \quad m+1\le k \le l \,,
\end{equation}
will be called type I summands.
If \eqref{108} is not of type I then let $k$, $m+1\le k \le l$,
be the minimal number such that
$$
\big(0 < s_k < \#\g_k\big)\quad \text{\rm\ or}\quad
\big(s_k=0\text{\rm\  and }s_{k+1}=\#\g_{k+1}\big) \,.
$$
We shall say that \eqref{108} is of type II (type III) if the 
first (second) parenthesis contains a true statement.

We shall show that the the type II summands cancel with
the type III summand while the type I summands produce
the RHS of \eqref{105}. 

The cancelation of the type II and type III summands 
follows from the following identity. Suppose
that
$$
0 < s_k < \#\g_k\,.
$$
Then
\begin{multline}\label{1010}
\sqb{\dots,\g_k,\dots}
{\dots,s_k,\dots} + \\
\sum_{\stackrel{\de \subset \g_k}{
\#\de=s_k}}
\left[\begin{matrix} \g_1,\dots,\g_m\\{}\end{matrix}
\left|\begin{matrix} 
\dots,&\g_k\setminus\de,&\de,&\dots \\ 
\dots,&0,&\#\de,&\dots
\end{matrix} \right. 
\right]
 = 0 \,. 
\end{multline}
To see \eqref{1010} notice that 
 all $\binom {\g_k}{s_k}$ summands in the sum over subsets
$\de\subset \g_k$ are equal and proportional
to the first summand in \eqref{1010}. 
 
Now consider a type I summand \eqref{109}. Set
\begin{equation}\label{1011}
\de=\g_{m+1}\cup \dots \cup \g_k \,. 
\end{equation}
We want to fix a subset 
$$
\de\subset \{2,\dots,n\}
$$
and compute the sum of 
all type I summands \eqref{109} satisfying \eqref{1011}.
Let us consider the nontrivial case 
$$
\de\ne\emptyset \,.
$$

First sum over all
$$
(\g_{m+1}, \dots, \g_k)\in\Ga(\de)\,,
$$ 
where $\Ga(\de)$ stands for the set of all compositions of
the set $\de$.
We have

\begin{multline}\label{1012}
\sqb
{\g_{m+1},&\dots,&\g_k,&\g_{k+1},\dots}
{\#\g_{m+1},&\dots,&\#\g_k,&0,\dots,0}=\\
(-1)^{k-m+1} \sqb
{\de &\g_{k+1},\dots}
{\#\de,&0,\dots,0} \,. 
\end{multline} 
Since  \eqref{1012} depends only on the parity of  number of parts 
in the composition \eqref{1011} we can use  \eqref{710}
to obtain 
\begin{multline*} 
\sum_{
(\g_{m+1},\dots,\g_k)\in\Ga(\de)}
\sqb
{\g_{m+1},&\dots,&\g_k,&\g_{k+1},\dots}
{\#\g_{m+1},&\dots,&\#\g_k,&0,\dots,0} = \\
(-1)^{\#\de+1}
\sqb
{\de, &\g_{k+1},\dots}
{\#\de,&0,\dots,0} \,.
\end{multline*}

Now sum over the remaining blocks 
$$
(\g_1,\dots,\g_{m},\g_{k+1},\dots,\g_l)\,.
$$
It follows from the definition \eqref{104}, \eqref{106} that
\begin{multline*}
(-1)^{\#\de+1}
\sum_{
(\g_1,\dots,\g_{m},\g_{k+1},\dots,\g_l)}
\sqb
{\de, &\g_{k+1},\dots}
{\#\de,&0,\dots,0}
 \\=
(-1)^{n+\ell(\si)} R^{\si} (t\ve t_0) \,,
\end{multline*}
where  $\si=\si(\de)$ is the following element of $\Po n$
$$
\si=\{{1\cup \de}, \text{Atom}(\{2,\dots,n\}\setminus \de)\} \,.
$$
Recall that  $\text{Atom}(\{2,\dots,n\}\setminus \de)$ denotes the 
partition of the set $\{2,\dots,n\}\setminus \de$ into
1-element blocks.
 
Thus,  the sum of all type I summands in \eqref{107}
equals the RHS of \eqref{105}. This  concludes the proof of 
the theorem. 
\end{proof}

\begin{cor}\label{cor102}
\begin{equation}\label{1014}
T(qt_1,t_2,\dots,t_n)= \sum_{\pi\in\Po n} 
(-1)^{n+\lpi} T^\pi(\tdn) \,. 
\end{equation}
\end{cor}
\begin{proof}
Take
$$
r(x;m)=\frac{\T^{(m)}(x)}{\T(x)} \,.
$$
By virtue of \eqref{62} this function satisfies \eqref{103} and 
it obviously satisfies \eqref{102}. We have
$$
T(\tdn)=\Res_{\,t_0=1} R(\tdn\ve t_0) \,.
$$
Taking the residue in \eqref{105} we obtain \eqref{1014}.

\end{proof}

\begin{cor}\label{cor103}
\begin{equation}\label{1015}
{U(qt_1,t_2,\dots,t_n)} = - q^{1/2}
t_1 \cdots t_n  \left( \sum_{\pi\in\Po n} 
(-1)^{n+\lpi} {U^\pi(\tdn)}
\right) \,. 
\end{equation}
\end{cor}

\section{Singularities of the RHS in \eqref{64}}\label{sec:Using}

In this section we shall prove that $U(\tdn)$
has exactly the same singularities as $F(\tdn)$.

\begin{thm}\label{thm111} For $k=1,\dots,n$ we have 
\begin{equation}\label{111}
U(t_1,\dots,t_n)=(-1)^{m}
 \frac{q^{m^2/2} m^{k-1}}{q^m t_1\cdots t_k-1}
\frac{U(t_{k+1},\dots,t_n)}{(t_{k+1}\cdots t_n)^m}+ \dots \,, 
\end{equation}
where dots stand for terms regular at the divisor 
$$
q^m t_1\cdots t_k=1 
$$ 
and we
assume that 
$$
U(t_{k+1},\dots,t_n)=1\,,\quad k=n\,.
$$
\end{thm}

Let us again point out that \eqref{111} implies that the 
residue
$$
\Res_{q^m t_1\dots t_k=1} U(\tdn)
$$
is independent of $t_1,\cdots,t_k$.
In the proof we shall use a curious identity good for
any odd smooth function which we shall state as a
separate Theorem \ref{thm113}.

It is clear that since the function $U(t)$ satisfies the
very same difference equation as $F(t)$ does it
suffices to consider the case
$$
m=0 \,.
$$
That is, we have to show that 
\begin{equation}\label{112}
U(\tdn)=\frac1{t_1-1} U(t_2,\dots,t_n)+ \dots\,,
\quad t_1\to 1 \,, 
\end{equation}
and that
\begin{equation}\label{113}
\text{\rm 
$U(t)$ is regular on $\{t_1\cdots t_k=1\}$ for
$1<k\le n$ 
} 
\end{equation}

It follows from the very definition of the 
function $U(t)$ (look, for example, at the
formula \eqref{66}) that
\begin{align*}
U(\tdn)&=U(t_1,\dots,t_k) U(t_{k+1}\dots t_n) + \dots \\
&=\frac{T(t_1,\dots,t_k)}{t_1 \cdots t_k -1}
U(t_{k+1}\dots t_n) + \dots\,,
\end{align*}
where dots stand for terms regular on the divisor 
$$
t_1\cdots t_k=1 \,.
$$ 
Since 
$$
T(1)=1 \,.
$$
the Theorem 4.1 will follow from the following

\begin{thm}\label{thm112} We have
\begin{equation}\label{114}
T(\tdn)\Big|_{t_1\cdots t_n=1}=0 
\end{equation}
provided $n>1$.

\end{thm}

\begin{proof}
Till the end of the proof the variables $\tdn$ will be 
always subject to constraint
\begin{equation}\label{115}
t_1\cdots t_n=1\,. 
\end{equation}

Induct on $n$. Since
$$
T(t_1,t_2)=
\frac{\T'(t_1)}{\T(t_1)} + \frac{\T'(t_2)}{\T(t_2)}
$$
the case $n=2$ is clear. 

First, show that
\begin{equation}\label{116}
\text{\rm $T(\tdn)$ is a constant.} 
\end{equation}
By the difference equation \eqref{1014} and the induction hypothesis
$$
T(q t_1, t_2,\dots,t_n)=T(\tdn)+(-1)^{n+1}  \,.
$$
Similarly
$$
T(q t_1, t_2,\dots,t_n)=T(q t_1, q^{-1}t_2,\dots,t_n)+(-1)^{n+1}  \,.
$$
It follows that
$$
T(q t_1, q^{-1}t_2,\dots,t_n)=T(\tdn) \,.
$$
Hence, by virtue of the strategy enunciated at the beginning of section \ref{sec:not}, the claim
\eqref{116} will follow from the following claim
\begin{equation}\label{117}
\text{\rm $T(\tdn)$ is regular}\,, 
\end{equation}
which will now be established. 

By symmetry it suffices to show that $T(\tdn)$ is
regular on the divisor
\begin{equation}\label{118}
t_1 \dots t_k = 1\,, \quad 1\le k < n-2 \,. 
\end{equation}
Again, from the definition \eqref{78} it is clear that
\begin{multline}\label{119}
T(\tdn)=\\
\left(\frac1{\T(t_1\dots t_k)} + 
\frac1{\T(t_{k+1}\dots t_n)} 
\right)T(t_1\dots t_k)T(t_{k+1}\dots t_n) + \dots 
\end{multline}
where dots stand for a function regular on \eqref{118}.
By \eqref{115} the sum in parentheses in \eqref{119} vanishes
and this proves \eqref{116} and \eqref{117}.  

Thus, what is left is to show that LHS of \eqref{114}
vanishes at some point. This follows from the
general identity established below in Theorem \ref{thm113} where one
has to
substitute
$$
f(x)=\T(x) \,.
$$
This concludes the proof of the theorem. 
\end{proof}

Let $f(x)$ be an odd function
$$
f(x^{-1})=-f(x) \,.
$$
Set
$$
\p{k}(x):= \left(x\d x\right)^k f(x)\,, k\in\N \,.
$$
Consider the following function 
\begin{equation}\label{1110}
\Phi(\tdn):=\sum_{\g\in\G n} \Psi_\g (t) \,, 
\end{equation}
where
$$
\Psi_\g(t)=(-1)^{\lga} \p{\#\g_1}(1)\, \prod_{i=2}^{\lga}
\frac{
\p{\#\g_i}\left(\tgi\right)}
{f\left(\tgi\right)} \,. 
$$
Note that $\Psi_\g(t)=0$ if $\#\g_1$ is even. We want to prove
the following

\begin{thm}\label{thm113} Suppose that $n>1$ and the variables $\tdn$ are subject
to constraint
$$
t_1 \cdots t_n = 1 \,. 
$$
Then for  any odd function $f(x)$ with a simple zero at $x=1$
$$
f(1)=0\,, \quad f'(1)\ne 0
$$
we have
\begin{equation}\label{1111}
\Phi(\tdn)\to 0\,, \quad t_1 \to 1 \,, 
\end{equation}
where the function $\Phi(t)$ was defined in \eqref{1110}.
\end{thm}

\begin{proof}
Let us divide all summands in \eqref{1110} into 3 following types: 
\begin{enumerate}
\item $\g_{\lga}=\{1\}$\,,
\item $\g_1=\{1\}$\,,
\item others.
\end{enumerate}
Note that all type 3 summands  are regular at $t_1=1$.

Let $\g=(\g_1,\dots,\g_{\ell-1},\{1\})$ be a type 1 composition.
Then the composition
$$
\g'=(\{1\},\g_1,\dots,\g_{\ell-1})
$$
is of type 2. Consider the sum
\begin{equation}\label{1112}
\Psi_\g (t) + \Psi_{\g'} (t) \,. 
\end{equation}
There are two possible cases: $\#\g_1$ is even and $\#\g_1$ is odd.
If
$$
\text{\rm $\#\g_1$ is even}
$$
then the first summand in \eqref{1112} is zero
and the second one is regular at $t_1=1$. In this case
we obtain
\begin{equation}\label{1114}
\lim_{t_1\to 1} \left(\Psi_\g (t) + \Psi_{\g'} (t)
\right) =
(-1)^{\ell} \p{\#\g_1+1}(1)\,\, \prod_{i=2}^{\ell-1}
\frac{
\p{\#\g_i}\left(\tgi\right)}
{f\left(\tgi\right)} \,. 
\end{equation}
It is easy to see that \eqref{1114} will exactly cancel with
the type 3 summand corresponding to the 
composition
\begin{equation}\label{1115}
(\{1\}\cup\g_1,\dots,\g_{\ell-1}) \,. 
\end{equation}
Note that if $\#\g_1$ is odd then the contribution of
the composition \eqref{1115} is zero. 

Now consider the sum \eqref{1112} in the case 
$$
\text{\rm $\#\g_1$ is odd}\,.
$$
We have
\begin{alignat*}{2}
\Psi_\g(t)&=
 (-1)^{\ell+1} \p{\#\g_1}(1)\,\frac{\pp(t_1)}{f(t_1)}
&&\prod_{i=2}^{\ell-1}
\frac{
\p{\#\g_i}\left(\tgi\right)}
{f\left(\tgi\right)}
\\
\Psi_{\g'}(t)&=
(-1)^{\ell} \pp(1) \,\frac{\p{\#\g_1}(t_1)}{f(t_1)}
 &&\prod_{i=2}^{\ell-1}
\frac{
\p{\#\g_i}\left(t_1\cdot \tgi\right)}
{f\left(t_1 \cdot\tgi\right)}
\end{alignat*}
Observe that 
\begin{equation}\label{1117}
\pp(1) \,\, \frac{\p{\#\g_1}(t_1)}{f(t_1)} -
\p{\#\g_1}(1)\,\, \frac{\pp(t_1)}{f(t_1)} \to 0\,,
\quad t_1 \to 1 
\end{equation}
because \eqref{1117} is regular and odd. Therefore
\begin{multline*}
\lim_{t_1\to 1} \left(\Psi_\g (t) + \Psi_{\g'} (t)
\right) = \\
(-1)^\ell
\p{\#\g_1}(1)\,\, 
\lim_{t_1\to 1}  \frac{\pp(t_1)}{f(t_1)}
\left(
\prod_{i=2}^{\ell-1}
\frac{
\p{\#\g_i}\left(t_1\cdot \tgi\right)}
{f\left(t_1 \cdot\tgi\right)} - \right.\\
\left.
\prod_{i=2}^{\ell-1}
\frac{
\p{\#\g_i}\left(\tgi\right)}
{f\left(\tgi\right)}\right) \,.
\end{multline*}
By L'Hospital rule this limit equals
\begin{multline*}
(-1)^\ell \p{\#\g_1}(1)
\left(
\sum_{k=2}^{\ell-1} 
\left( \frac
{\p{\#\g_k+1}\left({\tprod}_{j\in\ups(\g)_k} t_j\right)}
{f\left({\tprod}_{j\in\ups(\g)_k} t_j\right)} -
\right.
\right. \\
-
\left.
\left.
\frac
{\p{\#\g_k}\left({\tprod}_{j\in\ups(\g)_k} t_j\right)\pp\left({\tprod}_{j\in\ups(\g)_k} t_j\right)}
{f\left({\tprod}_{j\in\ups(\g)_k} t_j\right)^2}
\right) 
\prod_{\stackrel{i=2}{i\ne k}}^{\ell-1}
\frac{
\p{\#\g_i}\left(\tgi\right)}
{f\left(\tgi\right)}
\right) \,.
\end{multline*}
It is easy to see that the summand
$$
(-1)^\ell \p{\#\g_1}(1) \,\,
\frac
{\p{\#\g_k+1}\left({\tprod}_{j\in\ups(\g)_k} t_j\right)}
{f\left({\tprod}_{j\in\ups(\g)_k} t_j\right)} 
\prod_{\stackrel{i=2}{i\ne k}}^{\ell-1}
\frac{
\p{\#\g_i}\left(\tgi\right)}
{f\left(\tgi\right)}
$$
cancels with the contribution of the type 3
composition
\begin{equation}\label{1118}
(\g_1,\dots, \{1\}\cup \g_k ,\dots, \g_{\ell-1})
\,, \quad 2\le k \le l-1\,. 
\end{equation}
Similarly the summand 
\begin{multline*}
(-1)^{\ell+1} \p{\#\g_1}(1) \,\,
\frac
{\p{\#\g_k}\left({\tprod}_{j\in\ups(\g)_k} t_j\right)\pp\left({\tprod}_{j\in\ups(\g)_k} t_j\right)}
{f\left({\tprod}_{j\in\ups(\g)_k} t_j\right)^2} \times \\
\prod_{\stackrel{i=2}{i\ne k}}^{\ell-1}
\frac{
\p{\#\g_i}\left(\tgi\right)}
{f\left(\tgi\right)}
\,.
\end{multline*}
cancels with the contribution of the type 3
composition
\begin{equation}\label{1119}
(\g_1,\dots,\{1\},\g_k,\dots,\g_{\ell-1})
\,, \quad 2\le k \le l-1 \,. 
\end{equation}

It is clear that the compositions of the form
\eqref{1115}, \eqref{1118}, and \eqref{1119} exhaust the set of
type 3 compositions. This concludes the proof.
\end{proof}

\section{Conclusion of the proof of Theorem \ref{thm61}}\label{sec:pf}

Induct on $n$. Suppose that $n=1$. By Theorems \ref{thm81} and \ref{thm91}
the function
$$
\T(t_1) F(t_1)
$$
is holomorphic on $\C\setminus 0$, invariant
under the transformation
$$
t_1 \mapsto q t_1
$$
and equal to $1$ for $t_1=1$. It follows that
$$
F(t_1)=\frac{1}{\T(t_1)} \,.
$$

Suppose that $n>1$ and consider the function
\begin{equation}\label{121}
\T(t_1\cdots t_n) F(\tdn) - T(\tdn) \,. 
\end{equation}
By induction hypothesis, we have
$$
\T(t_1\cdots t_n) F^\pi (\tdn) - T^\pi(\tdn)\,,
$$
provided
$$
\lpi < n \,.
$$
Therefore by Theorem 1.1 and  Corollary 3.2  the function  
is invariant
under the transformation 
$$
t_1 \mapsto q t_1\,. 
$$
By symmetry it is invariant under all 
transformations 
$$
t_i \mapsto q t_i\,, \quad i=1,\dots,n \,. 
$$
By Theorems \ref{thm91} and \ref{thm111}  the function \eqref{121} is holomorphic and vanishes
if
$$
t_1 \cdots t_n = 0 \,.
$$
Therefore \eqref{121} equals zero. This concludes the proof.

\setcounter{section}{12}

\section{Another Example}\label{sec:example}

Another representation of a subalgebra of the algebra $
\sD$ of differential operators on
$\C[t,t^{-1}]$ consisting of differential operators which are skew-adjoint in a suitable sense,
was studied in \cite{B}. This algebra contains Virasoro and also the odd powers
$D,D^3,D^5,\ldots$ of $D = t\frac{d}{dt}$. The $D^{2n+1}$ act semisimple with finite eigenspaces,
and the resulting character is
\begin{equation}\label{131}\Psi(\tau_1,\tau_3,\ldots) = q_1^{\zeta(-1)/2}q_3^{\zeta(-3)/2}\cdots
\prod_{n=1}^\infty (1-q_1^nq_3^{n^3}\cdots)^{-1}
\end{equation}
In this section, we will show that the Taylor expansion for $\Psi$ is quasimodular of weight -1/2
(cf. \eqref{4wt}), and we will calculate the $n$-point function 
\begin{multline}\label{132}\sF_n(\tau_1,z_1,\dotsc,z_n) := \\
\sum_{k_1,\dotsc,k_n\ge 1}
\frac{\partial^n}{\partial\tau_{2k_1-1}\cdots\partial\tau_{2k_n-1}}(\Psi)|_{\tau_3 = \cdots = 0}
\frac{z_1^{2k_1-1}\cdots z_n^{2k_n-1}}{(2k_1-1)!\cdots (2k_n-1)!}.
\end{multline}
As before, we write 
$$q_r = \exp(2\pi i \tau_r);\quad r\ge 1.
$$

We compute
\begin{gather*} \frac{1}{2\pi i}\frac{\partial}{\partial \tau_{2j-1}} \log\Psi = \zeta(1-2j)/2 +
\sum_{m,n=1}^\infty m^{2j-1}q_1^{nm}q_3^{nm^3}\cdots \\
\frac{1}{(2\pi i)^r}\frac{\partial^r}{\partial \tau_{2j_1-1}\cdots\partial\tau_{2j_r-1}} \log\Psi
= 
\sum_{m,n=1}^\infty m^{2(\sum_k j_k-r)-1}(nm)^{r-1}q_1^{nm}q_3^{nm^3}\cdots \notag
\end{gather*}
Note these expressions depend only on $r$ and $\sum_{k=1}^r j_k$. Define
\begin{gather}\label{134} h_{r, 2(j_1+\ldots+j_r)}(\tau_1,\tau_3,\ldots) := \frac{1}{(2\pi
i)^r}\frac{\partial^r}{\partial \tau_{2j_1-1}\cdots\partial\tau_{2j_r-1}} \log\Psi,\\
g_{r, 2(j_1+\ldots+j_r)}(\tau_1) := h_{r, 2(j_1+\ldots+j_r)}(\tau_1,0,0,\ldots) =
\frac{\partial^{r-1}}{\partial\tau_1^{r-1}}G_{2(j_1+\ldots+j_r-r+1)}(\tau_1) \label{135}
\end{gather}
where $G_{2p}(\tau)$ is the Eisenstein series of weight $2p$ as in \eqref{34}. Note $g_{r,s}$ is
quasimodular of weight $2a$. For $A = (a_3,a_5,\ldots)$ write $\wt(A) = 4a_3+6a_5+\ldots$
(compare \eqref{4wt}). The coefficient of $\tau^A/A!$ in the Taylor expansion for $\log\Psi$ is 
$g_{r,\wt A}(\tau_1)$. Also
\begin{equation}\label{136} \Psi|_{\tau_3=\cdots = 0} = \eta^{-1}(\tau_1)
\end{equation}
Thus, by \eqref{4wt} we conclude that $\log(\Psi\eta)$ is quasimodular of
weight $0$. Exponentiating yields
\begin{prop}\label{prop131} $\Psi$ is quasimodular of weight $-1/2$. 
\end{prop} 

We now consider the $n$-point function \eqref{132}
\begin{lem} Let $\sS = \{j_1,\dotsc,j_n\}$ be a set of positive integers. Let $\Pi(\sS)$ denote
the set of all partitions $\mu = \{\mu_1,\dotsc,\mu_\ell\}$ of $\sS$. For any finite set $\phi$
let $\#\phi$ denote  the number of elements in $\phi$. For a subset $\phi\subset \sS$ let $|\phi|$
be the sum over the elements. Then
$$\frac{1}{(2\pi
i)^n}\frac{\partial^n}{\partial \tau_{2j_1-1}\cdots\partial\tau_{2j_n-1}}\Psi = \Psi\sum_{\mu\in
\Pi(\sS)} \prod_{k=1}^{\#\mu}h_{\#\mu_k,2|\mu_k|}
$$ 
\end{lem}
\begin{proof}By induction on $n$. For $n=1$ this is just the logarithmic derivative. Suppose now
$n\ge 2$ and the assertion holds for $n-1$. Note
$$\frac{1}{2\pi i}\frac{1}{\partial\tau_{2j_n-1}}h_{n-1,2(j_1+\ldots+j_{n-1})} =
h_{n,2(j_1+\ldots+j_{n})}
$$
Let $\sT = \{j_1,\dotsc,j_{n-1}\}$. We have inductively
\begin{multline*}\frac{1}{(2\pi
i)^n}\frac{\partial^n}{\partial \tau_{2j_1-1}\cdots\partial\tau_{2j_n-1}}\Psi =
\frac{1}{2\pi i}\frac{\partial}{\partial\tau_{2j_n-1}}\Big(\Psi\sum_{\mu\in
\Pi(\sT)} \prod_{k=1}^{\#\mu}h_{\#\mu_k,2|\mu_k|}\Big) = \\
= \Psi h_{1,2j_n}\sum_{\mu\in
\Pi(\sT)} \prod_{k=1}^{\#\mu}h_{\#\mu_k,2|\mu_k|} + \Psi\sum_{\mu\in
\Pi(\sT)}\sum_{p=1}^{\#\mu}h_{\#\mu_p+1,2(|\mu_p|+j_n)}\prod_{k\neq p}h_{\#\mu_k,2|\mu_k|} \\
= \Psi\sum_{\mu\in\Pi(\sS)}\prod_{k=1}^{\#\mu}h_{\#\mu_k,2|\mu_k|}.
\end{multline*}
\end{proof}

Define
\begin{equation}\sG_n(\tau,z) := \sum_{r=1}^\infty \frac{\partial^{n-1}}{\partial\tau^{n-1}}
G_{2r}(\tau)z^{2r+n-2}/(2r+n-2)!
\end{equation}
One has, with $\Theta$ as in \eqref{61}
\begin{equation}2\sG_1(\tau,2\pi iz) := \frac{1}{2\pi
i}(-\frac{d}{dz}\log\Theta(z)+\frac{1}{z}) 
\end{equation}
To see this, let $\sigma(z,\tau)$ be the elliptic sigma function as defined e.g. in \cite{L}, p.
247. One has
\begin{multline*}\frac{-1}{2\pi i}\frac{d}{dz}\log\Theta = \frac{-1}{2\pi i} \frac{d}{dz}\log\sigma
+2G_2(\tau)(2\pi iz) = \\
\frac{-1}{2\pi iz} +2G_2(\tau)(2\pi iz)+2G_4(\tau)(2\pi iz)^3/3!+2G_6(\tau)(2\pi iz)^5/5!+\ldots
\end{multline*}
Combining \eqref{132}, \eqref{135}, and \eqref{136} we get
\begin{multline}\label{1310} \sF(\tau,z_1,\dotsc,z_n) = \eta(\tau)^{-1}\sum_{k_1,\ldots,k_n\ge 1}
\frac{z_1^{2k_1-1}\cdots z_n^{2k_n-1}}{(2k_1-1)!\cdots (2k_n-1)!}\times \\
\sum_{\mu\in\Pi\{k_1,\dotsc,k_n\}} \prod_{p=1}^{\#\mu} \frac{\partial^{\#\mu_p
-1}}{\partial\tau^{\#\mu_p -1}}G_{2(|\mu_p|-\#\mu_p +1)}(\tau).
\end{multline}

For a function $f(z_1,\dotsc,z_n)$ let $\epsilon f$ denote the ``oddification'' of $f$, e.g.
$$\epsilon f(z_1,z_2) = \frac{1}{4}(f(z_1,z_2)-f(-z_1,z_2)-f(z_1,-z_2)+f(-z_1,-z_2))
$$
Note
$$\epsilon((z_1+\ldots+z_n)^r/r!) = \sum_{\ell_1+\ldots+\ell_n=(r+n)/2}
\frac{z_1^{2\ell_1-1}\cdots z_n^{2\ell_n-1}}{(2\ell_1-1)!\cdots (2\ell_n-1)!}
$$
It follows that
\begin{multline*}\epsilon\sG_n(\tau,z_1+\ldots+z_n) = \sum_{r=1}^\infty \frac{\partial^{n-1}}
{\partial\tau^{n-1}} G_{2r}(\tau)\times \\
\sum_{k_1+\ldots+k_n=r+n-1}\frac{z_1^{2k_1-1}\cdots
z_n^{2k_n-1}}{(2k_1-1)!\cdots (2k_n-1)!}
\end{multline*}
We conclude
\begin{prop} The $n$-point function is given by
$$\sF(\tau,z_1,\dotsc,z_n) =
\eta(\tau)^{-1}\sum_{\mu\in\Pi_n}\prod_{k=1}^{k=\#\mu}\epsilon\sG_{\#\mu_k}(\sum_{i\in \mu_k} z_i)
$$
\end{prop}
\begin{proof} We have an obvious identification
$$\Pi_n = \Pi\{1,\dotsc,n\} = \Pi\{k_1,\dotsc,k_n\}
$$
(Note $\{k_1,\dotsc,k_n\}$ is treated formally, i.e. we do not identify the $k_i$ if they are
equal.) Write $\kappa=\{k_1,\dotsc,k_n\}$ and let $\mu_{p,\kappa}\subset \kappa$ be the
corresponding subset. Given $\mu\in\Pi_n$ define
\begin{multline*}\sF_\mu(\tau,z_1,\dotsc,z_n) = \eta(\tau)^{-1}\sum_{k_1,\ldots,k_n\ge 1}
\frac{z_1^{2k_1-1}\cdots z_n^{2k_n-1}}{(2k_1-1)!\cdots (2k_n-1)!}\times \\
 \prod_{p=1}^{\#\mu} \frac{\partial^{\#\mu_p
-1}}{\partial\tau^{\#\mu_p -1}}G_{2(|\mu_{p,\kappa}|-\#\mu_p +1)}(\tau).
\end{multline*} 
so $\sF = \sum_{\Pi_n} \sF_\mu$. It will suffice to show
\begin{equation}\label{1312}\sF_\mu =
\eta(\tau)^{-1}\prod_{p=1}^{p=\#\mu}\epsilon\sG_{\#\mu_p}(\sum_{i\in \mu_p} z_i).
\end{equation}
But given $\kappa$, the term $z_1^{2k_1-1}\cdots z_n^{2k_n-1}$ occurs exactly once in
\eqref{1312} and has the correct coefficient in $\tau$. 
\end{proof}

\bibliographystyle{plain}
\renewcommand\refname{References}

\end{document}